\documentclass[12pt,leqno]{article}
\usepackage{amssymb}
\usepackage[mathscr]{eucal}
\usepackage{amsmath,amssymb,latexsym,theorem,enumerate,ifthen}
\usepackage[english]{babel}
\usepackage{color,url}
\usepackage{appendix}
\usepackage{graphicx}
\usepackage{subcaption}

\newcommand{\NN}{\mathbb{N}}
\newcommand{\QQ}{\mathbb{Q}}
\newcommand{\RR}{\mathbb{R}}

\newcommand{\ZZ}{\mathbb{Z}}

\newcommand{\bx}{{\boldsymbol{x}}}

\newcommand{\blambda}{{\boldsymbol{\lambda}}}

\newcommand{\bmu}{{\boldsymbol{\mu}}}

\newcommand{\bone}{{\boldsymbol{1}}}

\newcommand{\cA}{{\mathcal A}}
\newcommand{\cB}{{\mathcal B}}

\newcommand{\cE}{{\mathcal E}}

\newcommand{\cH}{{\mathcal H}}

\newcommand{\cX}{{\mathcal X}}

\newcommand{\dd}{\mathrm{d}}
\newcommand{\ee}{\mathrm{e}}

\DeclareMathOperator*{\conv}{conv}

\newcommand{\EE}{\operatorname{\mathbb{E}}}

\newcommand{\PP}{\operatorname{\mathbb{P}}}

\newcommand{\Med}{\operatorname{Med}}

\newcommand{\sign}{\operatorname{sign}}

\newcommand{\vare}{\varepsilon}
\newcommand{\comment}[1]{}

\renewcommand{\mid}{\,|\,}

\renewcommand{\leq}{\leqslant}
\renewcommand{\geq}{\geqslant}

\newcommand{\proofend}{\hfill\mbox{$\Box$}}
\newcommand{\Mid}{\,\bigg|\,}

\numberwithin{equation}{section}

\theoremstyle{change} \theorembodyfont{\em}
\newtheorem{Lem}{Lemma.}[section]
\newtheorem{Thm}[Lem]{Theorem.}
\newtheorem{Pro}[Lem]{Proposition.}

\newtheorem{Cor}[Lem]{Corollary.}
\newtheorem{Def}[Lem]{Definition.}

\theorembodyfont{\rm}
\newtheorem{Rem}[Lem]{Remark.}

\newtheorem{Ex}[Lem]{Example.}

\def\OnlyOnArXiv#1#2{\ifthenelse{\equal{#1}{Y}}{#2}{}}

\newenvironment{proof}{\noindent{\bf Proof.}}{\proofend}

\begin{document}

\begin{center}
 {\bfseries\Large Existence and uniqueness\\[1mm] of weighted generalized $\psi$-estimators}

\vspace*{3mm}

{\sc\large
  M\'aty\'as $\text{Barczy}^{*}$,
  Zsolt $\text{P\'ales}^{**,\diamond}$ }

\end{center}

\vskip0.2cm

\noindent
 * HUN-REN–SZTE Analysis and Applications Research Group,
   Bolyai Institute, University of Szeged,
   Aradi v\'ertan\'uk tere 1, H--6720 Szeged, Hungary.

\noindent
 ** Institute of Mathematics, University of Debrecen,
    Pf.~400, H--4002 Debrecen, Hungary.

\noindent E-mails: barczy@math.u-szeged.hu (M. Barczy),
                  pales@science.unideb.hu  (Zs. P\'ales).

\noindent $\diamond$ Corresponding author.

\vskip0.2cm


{\renewcommand{\thefootnote}{}
\vspace*{0.2cm}
\footnote{M\'aty\'as Barczy was supported by the project TKP2021-NVA-09. Project no.\ TKP2021-NVA-09 has been implemented with the support provided by the Ministry of Culture and Innovation of Hungary from the National Research, Development and Innovation Fund, financed under the TKP2021-NVA funding scheme.
Zsolt P\'ales is supported by the K-134191 NKFIH Grant.}}

\vspace*{-10mm}

\begin{abstract}
We introduce the notions of generalized and weighted generalized $\psi$-estimators as unique points of sign change of some appropriate functions, and we give necessary as well as sufficient conditions for their existence.
We also derive a set of sufficient conditions under which the so-called $\psi$-expectation function has a unique point of sign change.
We present several examples from statistical estimation theory, where our results are well-applicable.
For example, we consider the cases of empirical quantiles, empirical expectiles, some $\psi$-estimators that are important in robust statistics, and some examples from maximum likelihood theory as well.
Further, we introduce Bajraktarevi\'c-type (in particular, quasi-arithmetic-type) $\psi$-estimators.
Our results specialized to $\psi$-estimators with a function $\psi$ being continuous in its second variable provide new results for (usual) $\psi$-estimators (also called Z-estimators).
\end{abstract}


\section{Introduction}
\label{Section_intro}

In statistics, M-estimators play a fundamental role, and a special subclass, the class of $\psi$-estimators (also called $Z$-estimators), is also in the heart of investigations.
The M-estimators (where the letter M refers to ''maximum likelihood-type'') were introduced by Huber \cite{Hub64, Hub67}.
Let $(X,\cX)$ be a measurable space, $\Theta$ be a Borel subset of $\RR$, and $\varrho:X\times\Theta\to\RR$
 be a function such that for each $t\in\Theta$, the function $X\ni x\mapsto \varrho(x,t)$ is measurable
 with respect to the sigma-algebra $\cX$.
Let $(\xi_k)_{k\geq 1}$ be a sequence of i.i.d.\ random variables with values in $X$ such that the distribution of $\xi_1$ depends on an unknown parameter $\vartheta \in\Theta$.
For each $n\geq 1$, Huber \cite{Hub64, Hub67} introduced an estimator of $\vartheta$ based on the observations
$\xi_1,\ldots,\xi_n$ as a solution $\widehat\vartheta_n:=\widehat\vartheta_n(\xi_1,\ldots,\xi_n)$ of the following minimization problem:
 \begin{align}\label{help_M_est_min_problem}
   \inf_{t\in\Theta}\sum_{i=1}^n \varrho(\xi_i,t),
 \end{align}
provided that such a solution exists.
One calls $\widehat\vartheta_n$ an M-estimator of the unknown parameter $\vartheta\in\Theta$ based on
 the i.i.d.\ observations $\xi_1,\ldots,\xi_n$.
In stochastic optimization, $\widehat\vartheta_n$ is called a Sample Average Approximation (SAA) of a solution of the classical risk neutral stochastic program $\inf_{t\in\Theta} \EE(\varrho(\xi_1,t))$, see, e.g., Shapiro et al. \cite[Chapter 5]{ShaDenRus}.
For historical fidelity, we note that Huber \cite{Hub64} considered the special case when $X:=\RR$, $\Theta:=\RR$, and the function $\varrho$ depends only on $x-t$,
 i.e., $\varrho(x,t):=f(x-t)$, $x\in\RR$, $t\in\Theta$,
 with some given nonconstant function $f:\RR\to\RR$.
Turning back to the general case, under suitable regularity assumptions, the minimization problem  \eqref{help_M_est_min_problem} can be solved by setting the derivative of the objective function (with respect to the unknown parameter) equal to zero:
 \[
 \sum_{i=1}^n \partial_2\varrho(\xi_i,t)=0, \qquad t\in\Theta,
 \]
where $\partial_2 \varrho$ denotes the (partial) derivative of $\varrho$ with respect to its second variable.
In the statistical literature, $\partial_2\varrho$ is often denoted by $\psi$, and hence in this case an M-estimator is often called a $\psi$-estimator, while other authors call it a Z-estimator (the letter Z refers to ''zero'').
For a detailed exposition of M-estimators and $\psi$-estimators (Z-estimators), see, e.g., Kosorok \cite[Sections 2.2.5 and 13]{Kos} or van der Vaart \cite[Section 5]{Vaa}.

Throughout this paper, let $\NN=\ZZ_{++}$, $\ZZ_+$, $\QQ$, $\RR$, $\RR_+$ and $\RR_{++}$ denote the sets of positive integers, non-negative integers, rational numbers, real numbers, non-negative
 real numbers and positive real numbers, respectively.
For a real number $y\in\RR$, its positive and negative parts as well as its upper and lower integer parts are denoted by $y^+$ and  $y^-$ as well as by $\lceil y\rceil$ and $\lfloor y\rfloor$, respectively.
For a subset $S\subseteq \RR$, the convex hull of $S$ (which is the smallest interval containing $S$) is denoted by $\conv(S)$.
For each $n\in\NN$, let us also introduce the set $\Lambda_n:=\RR_+^n\setminus\{(0,\ldots,0)\}$.
All the random variables are defined on an appropriate probability space $(\Omega,\cA,\PP)$.

To the best of our knowledge, the topic of existence and uniqueness of $\psi$-estimators is less addressed in the statistical literature.
In the present paper, we are going to investigate two basic problems of this field as presented below.
Roughly speaking, Problem 1 is about the existence and uniqueness of the newly introduced notions: generalized $\psi$-estimators and weighted generalized $\psi$-estimators.
Problem 2 is devoted to the existence and uniqueness of a point of sign change for so-called $\psi$-expectation functions.

{\sl Problem 1.} Let $X$ be a nonempty set, $\Theta$ be a nonempty open interval of $\RR$.
Let $\Psi(X,\Theta)$ denote the class of real-valued functions $\psi:X\times\Theta\to\RR$ such that,
 for each $x\in X$, there exist $t_+,t_-\in\Theta$ such that $t_+<t_-$ and $\psi(x,t_+)>0>\psi(x,t_-)$.
Roughly speaking, a function $\psi\in\Psi(X,\Theta)$ satisfies the following property:
 for each $x\in X$, the function $t\ni\Theta\mapsto \psi(x,t)$ changes sign (from positive to negative)
 on the interval $\Theta$ at least once.
Given a function $\psi\in\Psi(X,\Theta)$, $n\in\NN$ and $\pmb{x}=(x_1,\ldots,x_n)\in X^n$, let us consider the equation
 \begin{align}\label{psi_est_equation}
  \psi_{\pmb{x}}(t):=\sum_{i=1}^n \psi(x_i,t) = 0, \qquad t\in \Theta.
 \end{align}
More generally, for $n\in\NN$, $\pmb{x}=(x_1,\ldots,x_n)\in X^n$ and $\pmb{\lambda}=(\lambda_1,\ldots,\lambda_n)\in\Lambda_n$, we also consider the weighted equation
 \begin{align}\label{weighted_psi_estimator}
   \psi_{\pmb{x},\pmb{\lambda}}(t):=\sum_{i=1}^n \lambda_i \psi(x_i,t)=0,\qquad t\in \Theta.
 \end{align}
The basic question we are going to investigate now is to find necessary as well as sufficient conditions
for the unique solvability of the equations \eqref{psi_est_equation} and \eqref{weighted_psi_estimator}, respectively.
In a broader context, we are going to find necessary as well as sufficient conditions for the existence of a point of sign change
 (see Definition \ref{Def_sign_change}) for the functions $\psi_{\pmb{x}}$ and $\psi_{\pmb{x},\pmb{\lambda}}$
 introduced by \eqref{psi_est_equation} and \eqref{weighted_psi_estimator}, respectively.
It will turn out that the points of sign change in question are unique provided that they exist,
 and one can call them as a {\sl generalized $\psi$-estimator} and {\sl weighted generalized $\psi$-estimator},
 respectively, for some unknown parameter in $\Theta$ based
 on the realization $(x_1,\ldots,x_n)\in X^n$ and weights $(\lambda_1,\ldots,\lambda_n)\in\Lambda_n$.
In Proposition \ref{Pro_measurability}, we also study the measurability of (weighted) generalized $\psi$-estimators, provided that $X$ is a measurable space,
 but we emphasize that in our general setup, $X$ is not necessarily a measurable space,
 it can be an arbitrary nonempty set.
Concerning Problem 1, it does not matter whether the random variables $\xi_1,\ldots,\xi_n$ of which $(x_1,\ldots,x_n)$ is a realization are i.i.d.\ or not.
As future research, one could investigate the asymptotic properties of the (weighted) generalized $\psi$ estimators based on $(\xi_1,\ldots,\xi_n)$ as $n\to\infty$, when the property i.i.d.\ for the sequence $(\xi_k)_{k\geq 1}$ could play a role.

{\sl Problem 2.}
Let $(X,\cX)$ be a measurable space, $\Theta$ be a nonempty open interval of $\RR$, and $\psi:X\times \Theta\to\RR$ be a measurable function in its first variable, i.e., for each $t\in\Theta$, the mapping $X\ni x\mapsto \psi(x,t)$ is measurable with respect to the sigma-algebra $\cX$.
Further, let $\xi:\Omega\to X$ be a random variable defined on a probability space $(\Omega,\cA,\PP)$ such that $\EE(\vert\psi(\xi,t)\vert)<\infty$ for each $t\in\Theta$.
We investigate the question of existence of a unique point of sign change (see Definition \ref{Def_sign_change}) for the function
 \begin{align}\label{fgv_elmeleti}
    \Theta \ni t\mapsto  \EE(\psi(\xi,t)).
 \end{align}
In the literature, we could not find a name for the function \eqref{fgv_elmeleti},
 however, we may call it a $\psi$-expectation function.
Under appropriate conditions, the $\psi$-estimator (Z-estimator) based on i.i.d.\ observations $\xi_1,\ldots,\xi_n$ (where $\xi_1$ has the same law as that of $\xi$)
 is supposed to 'well-estimate' the zero of the function \eqref{fgv_elmeleti},
 provided that it exists uniquely, for more details, see, e.g., Kosorok \cite[Sections 2.2.5 and 13]{Kos}.

In what follows, we discuss the connections between the Problems 1 and 2 introduced above.
From a stochastic optimization point of view, Problem 1 can be considered as a Sample Average Approximation (SAA) of Problem 2 provided that $(x_1,\ldots,x_n)$ is a realization of i.i.d.\ random variables $\xi_1,\ldots,\xi_n$ (where $\xi_1$ has the same law as that of $\xi$).
Under appropriate conditions, the generalized $\psi$-estimator based on i.i.d.\ observations $\xi_1,\ldots,\xi_n$ (where $\xi_1$ has the same law as that of $\xi$) is supposed to 'well-estimate' the point of sign change of the function \eqref{fgv_elmeleti}, provided that it exists.
In this paper, we do not investigate this question, it could be a topic of future research.
Further, note that if $\xi$ is a simple random variable such that $\PP(\xi=x_i)=p_i$, $i=1,\ldots,n$, where $n\in\NN$,
 $(x_1,\ldots,x_n)\in X^n$ and $p_1,\ldots,p_n\geq 0$, $p_1+\cdots+p_n=1$,
 then $\EE(\psi(\xi,t)) = \sum_{i=1}^n p_i\psi(x_i,t)$, $t\in\Theta$, and hence, in this special case,
 Problem 2 is a special case of Problem 1.

To mention some papers related to Problem 1, we can refer, for example, to Huber \cite[Lemma 1]{Hub64}, Tibshirani \cite{Tib}
 and Ali and Tibshirani \cite{AliTib}.
Tibshirani \cite{Tib} considered the lasso (least absolute shrinkage and selection operator) problem, which is also known as the $\ell_1$-penalized linear regression.
The lasso estimator is a popular tool in the theory of  sparse linear regression, mathematically, it is a solution of a not necessarily strictly convex minimization problem, where a penalty term being the $\ell_1$-norm of the coefficient vector comes into play.
Tibshirani \cite{Tib} studied the question of uniqueness of the lasso estimator.
Recently, Ali and Tibshirani \cite{AliTib} have studied the uniqueness of a generalized lasso estimator, where the penalty term in the corresponding minimization problem is the $\ell_1$-norm of a (penalty) matrix times the coefficient vector.
To mention further papers related to Problem 2, we can refer to Huber \cite[Lemma 2]{Hub64}, Clarke \cite{Cla}
(for details, see Remark \ref{Rem_Clarke}), Mathieu \cite{Mat} (for details, see Examples \ref{Ex_Mat_elmeleti} and \ref{Ex_Mat}),
 and to the very recent paper of Dimitriadis et al.\ \cite[Propositions S1, S2 and S3]{DimFisZie}, in which the authors,
 in particular, considered solvability of the equation $\EE(\psi(\zeta,\eta,t))=0$, $t\in\Theta$,
 where $\psi:\RR\times\RR\times\Theta\to\RR$ is a measurable function, $\Theta$ is a (non-empty) open parameter set
 of $\RR$, and $(\zeta,\eta)$ is a response-regression  (covariate) pair in some regression model.

Section \ref{Section_Ex_uniq_empirical} is about the existence and uniqueness of (weighted) generalized $\psi$-estimators (Problem 1).
First, we introduce the required terminology: the notions of point of sign change and level of increase
 for a real-valued function defined on a nonempty open interval (see Definitions \ref{Def_sign_change} and \ref{Def_level_of_increase}),
 and the notions of properties $[T_n]$ and $[T_n^{\pmb{\lambda}}]$ for a function in $\Psi(X,\Theta)$ (see Definition \ref{Def_Tn}, but we also present below).
We say that a function $\psi\in\Psi(X,\Theta)$ \emph{possesses the property
 $[T_n^{\pmb{\lambda}}]$ for some $n\in\NN$ and $\pmb{\lambda}=(\lambda_1,\ldots,\lambda_n)\in\Lambda_n$}
 if there exists a mapping $\vartheta_{n,\psi}^{\pmb{\lambda}}:X^n\to\Theta$ such that,
 for each $\pmb{x}=(x_1,\dots,x_n)\in X^n$ and $t\in\Theta$,
 \begin{align*}
    \psi_{\pmb{x},\pmb{\lambda}}(t)= \sum_{i=1}^n \lambda_i\psi(x_i,t)
     \begin{cases}
          > 0 & \text{if $t<\vartheta_{n,\psi}^{\pmb{\lambda}}(\pmb{x})$,}\\
          < 0 & \text{if $t>\vartheta_{n,\psi}^{\pmb{\lambda}}(\pmb{x})$}.
     \end{cases}
 \end{align*}
Note that if there exists such a mapping $\vartheta_{n,\psi}^{\pmb{\lambda}}$, then it is unique.
In case of $\lambda_i=1$, $i=1,\ldots,n$ (or equivalently, in case of equal positive weights), the property $[T_n^{\pmb{\lambda}}]$ is called property $[T_n]$.
In the first main result of our paper (see Theorem \ref{Thm_M_est_uniq} below), necessary as well as sufficient conditions are given for the properties $[T_n]$ and $[T_n^{\pmb{\lambda}}]$.
If $\psi$ is continuous in its second variable as well, then such conditions imply the unique existence of the corresponding (usual) $\psi$-estimator.
After Theorem \ref{Thm_M_est_uniq}, we present some properties of the property $[T_n^{\pmb{\lambda}}]$.
For example, Proposition \ref{Pro_3} is about a connection between the property $[T_n]$ and the strictly $\frac1n$-increasingness of some appropriately defined functions,
 and in Proposition \ref{T_n} we establish a 'grouping' property of the property $[T_n^{\pmb{\lambda}}]$.
Examples \ref{Ex_T2_fail}, \ref{Ex_not_omitted} and \ref{Ex_div_noomit} highlight the role of the
 conditions in Theorem \ref{Thm_M_est_uniq} and in Proposition \ref{T_n}.
We also introduce a class of Bajraktarevi\'c-type (in particular, quasi-arithmetic-type) $\psi$-estimators (motivated by the representation of Bajraktarevi\'c means as special deviation means) for which our results are well-applicable, see Definition \ref{Ex_Bajrak_emp} and Proposition \ref{Pro_BajType_est}.

Section \ref{Section_Ex_uniq_theoretical} is devoted to study the existence and uniqueness of the point of sign change of the $\psi$-expectation function given in \eqref{fgv_elmeleti} (Problem 2).
As the second main result of our paper, in Theorem \ref{Thm_theor_M_est_uniq}, we give a set of sufficient conditions in order that the  $\psi$-expectation function in question have a unique point of sign change.
We apply our results for $\psi$-expectation functions, when $\psi$ is a Bajraktarevi\'c-type function (see Proposition \ref{Pro_Bajrak_theor}), and when $\psi$ has the form used by Mathieu \cite{Mat} (see Example \ref{Ex_Mat_elmeleti} and
Proposition \ref{Pro_Math_type}). In Remark \ref{Rem_Clarke}, restricted to a one-dimensional parameter set,
 we recall Theorem 3.2 in Clarke \cite{Cla} on the local uniqueness of a root of
 the function \eqref{fgv_elmeleti}.
We will see that the assumptions of Theorem 3.2 in Clarke \cite{Cla} are much more involved
 and quite different compared to those of our Theorem \ref{Thm_theor_M_est_uniq}.

In Section \ref{Section_statistical_examples} we present several examples from statistical estimation theory that demonstrate the applicability
 of our results in Sections \ref{Section_Ex_uniq_empirical} and \ref{Section_Ex_uniq_theoretical}.
These examples may be divided into three main groups.
The first group of examples includes several well-known descriptive
statistics that can be considered as special $\psi$-estimators.
Namely, the empirical median (Example \ref{Ex_median}), the empirical quantiles (Example \ref{Ex_quantile})
 and the empirical expectiles (Example \ref{Ex_expectile}).
In particular, in Proposition \ref{Pro_quantile} we show that, given $n\geq 2$, the function $\psi$ corresponding
 to the empirical $\alpha$-quantile has the property $[T_n]$ if and only if $\alpha\notin\{\frac{1}{n},\ldots,\frac{n-1}{n}\}$.
The second group of examples contains the class of $\psi$-estimators recently used by Mathieu \cite{Mat}
 (Example \ref{Ex_Mat}), and some $\psi$-estimators that are important in robust statistics.
In particular, in Proposition 4.4, we derive necessary and sufficient conditions under which
 the function $\psi$ used in Mathieu \cite{Mat} has the property $[T_n^{\blambda}]$ for each $n\in\NN$ and $\blambda\in\Lambda_n$.
We emphasize that in all the above examples, we investigate the existence and uniqueness of (weighted) generalized $\psi$-estimators,
 compared to the existing results that addressed $\psi$-estimators (Z-estimators).
The third group of examples demonstrates the applicability of Theorem \ref{Thm_M_est_uniq} together
 with Proposition \ref{Pro_M_est_uniq} for proving existence and uniqueness of solutions of likelihood equations.
In Example \ref{Ex_1_norm}, we consider the maximum likelihood estimator (MLE) of one of the parameters of a normally distributed random
 variable supposing that its other parameter is known.
In Example \ref{Ex_1_norm_kevert}, we consider a mixture density function of the standard normally density function and the density function of a normally distributed random variable with mean $m\in\RR$ and variance $\sigma^2>0$ with equal weights $\frac{1}{2}$, and we study the solutions of the likelihood equation for  $m$ provided that $\sigma$ is known.

Sections \ref{Section_Proof_for_Section2}, \ref{Section_Proof_for_Section3} and \ref{Section_Proof_for_Section4} are devoted to the proofs
 of the results in Sections \ref{Section_Ex_uniq_empirical}, \ref{Section_Ex_uniq_theoretical} and \ref{Section_statistical_examples},
 respectively.

Now, we summarize the novelties of the paper.
We extensively discuss that, up to our knowledge, only few results are available for the existence and uniqueness of $\psi$-estimators and of the roots of $\psi$-expectation functions, and our paper can be considered as a theoretical contribution to this field.
Another important feature of our paper is that we present a broad variety of examples from statistical estimation theory, where our results can be well-applied.

In the end, we mention that, in the literature, one can find $\psi$-estimators, where the function $\psi$ depends on the sample size $n$ as well, see, e.g., Hampel, Hennig and Ronchetti \cite[Section 2]{HamHenRon}.
As future research, one might try to generalize the notion of (weighted) generalized $\psi$-estimators and our results to this more general case.
Another possible direction for future research is to explore the extension of our setup and results from a one-dimensional parameter set $\Theta$ to a multidimensional one (note that, in our present setup, $\Theta$ is supposed to be a nonempty open interval of $\RR$).

\section{Notions and results on the existence and uniqueness of weighted generalized $\psi$-estimators}
\label{Section_Ex_uniq_empirical}

To investigate Problem 1 presented in the Introduction, we introduce the required terminology.
First, we introduce the notion of a point of sign change for real-valued functions defined on an open interval.

\begin{Def}\label{Def_sign_change}
Let $\Theta$ be a nonempty open interval of $\RR$. For a function $f:\Theta\to\RR$, consider the following three level sets
\[
  \Theta_{f>0}:=\{t\in \Theta: f(t)>0\},\qquad
  \Theta_{f=0}:=\{t\in \Theta: f(t)=0\},\qquad
  \Theta_{f<0}:=\{t\in \Theta: f(t)<0\}.
\]
We say that $\vartheta\in\Theta$ is a \emph{point of sign change (of decreasing type) for $f$} if
 \[
 f(t) > 0 \quad \text{for $t<\vartheta$,}
   \qquad \text{and} \qquad
    f(t)< 0 \quad  \text{for $t>\vartheta$.}
 \]
\end{Def}

\begin{Rem}\label{Rem_2}
Note that, if $\vartheta\in\Theta$ is a point of sign change for $f$, then $\Theta_{f>0}$ and $\Theta_{f<0}$ are nonempty sets and $\sup\Theta_{f>0}=\inf\Theta_{f<0}=\vartheta$. Furthermore, there can exist at most one element $\vartheta\in\Theta$ which is a point of sign change for $f$.
 If $f$ is continuous at a point $\vartheta$ of sign change, then $f(\vartheta)=0$, moreover $\Theta_{f=0}=\{\vartheta\}$.
 Conversely, as an easy consequence of the Bolzano theorem, if $f:\Theta\to\RR$ is continuous,
 the sets $\Theta_{f>0}$ and $\Theta_{f<0}$ are nonempty and $f$ has a unique zero $\vartheta\in\Theta$,
 i.e., $\Theta_{f=0}=\{\vartheta\}$ holds, then $\vartheta$ is a point of sign change either for $f$ or for $(-f)$.
The continuity of $f$, however, is not necessary for the existence of a point of sign change for $f$.
For example, if $f$ is strictly decreasing and the sets $\Theta_{f>0}$ and $\Theta_{f<0}$ are nonempty,
 then it is easy to see that there exists a point of sign change for $f$.
\proofend
\end{Rem}

\begin{Def}\label{Def_Tn}
We say that a function $\psi\in\Psi(X,\Theta)$
  \begin{enumerate}[(i)]
    \item \emph{possesses the property $[C]$ (briefly, $\psi$ is a $C$-function)} if
           it is continuous in its second variable, i.e., if, for all $x\in X$,
           the mapping $\Theta\ni t\mapsto \psi(x,t)$ is continuous.
    \item \emph{possesses the property $[T_n]$ (briefly, $\psi$ is a $T_n$-function)
           for some $n\in\NN$} if there exists a mapping $\vartheta_{n,\psi}:X^n\to\Theta$ such that,
           for all $\pmb{x}=(x_1,\dots,x_n)\in X^n$ and $t\in\Theta$,
           \begin{align*}
             \psi_{\pmb{x}}(t):=\sum_{i=1}^n \psi(x_i,t) \begin{cases}
                 > 0 & \text{if $t<\vartheta_{n,\psi}(\pmb{x})$,}\\
                 < 0 & \text{if $t>\vartheta_{n,\psi}(\pmb{x})$},
            \end{cases}
           \end{align*}
          that is, for all $\pmb{x}\in X^n$, the value $\vartheta_{n,\psi}(\pmb{x})$ is a point of sign change for the function $\psi_{\pmb{x}}$. If there is no confusion, instead of $\vartheta_{n,\psi}$ we simply write $\vartheta_n$.
          We may call $\vartheta_{n,\psi}(\pmb{x})$ as a generalized $\psi$-estimator for
         some unknown parameter in $\Theta$ based on the realization $\bx=(x_1,\ldots,x_n)\in X^n$. If, for each $n\in\NN$, $\psi$ is a $T_n$-function, then we say that \emph{$\psi$ possesses the property $[T]$ (briefly, $\psi$ is a $T$-function)}.
    \item \emph{possesses the property $[Z_n]$ (briefly, $\psi$ is a $Z_n$-function) for some $n\in\NN$} if it is a $T_n$-function and
    \[
   \psi_{\pmb{x}}(\vartheta_{n,\psi}(\pmb{x}))=\sum_{i=1}^n \psi(x_i,\vartheta_{n,\psi}(\pmb{x}))= 0
    \qquad \text{for all}\quad \pmb{x}=(x_1,\ldots,x_n)\in X^n.
    \]
    If, for each $n\in\NN$, $\psi$ is a $Z_n$-function, then we say that \emph{$\psi$ possesses the property $[Z]$ (briefly, $\psi$ is a $Z$-function)}.
    \item \emph{possesses the property $[T_n^{\pmb{\lambda}}]$ for some $n\in\NN$ and $\pmb{\lambda}=(\lambda_1,\ldots,\lambda_n)\in\Lambda_n$ (briefly, $\psi$ is a $T_n^{\pmb{\lambda}}$-function)} if there exists a mapping $\vartheta_{n,\psi}^{\pmb{\lambda}}:X^n\to\Theta$ such that, for all $\pmb{x}=(x_1,\dots,x_n)\in X^n$ and $t\in\Theta$,
          \begin{align*}
           \psi_{\pmb{x},\pmb{\lambda}}(t):= \sum_{i=1}^n \lambda_i\psi(x_i,t) \begin{cases}
                 > 0 & \text{if $t<\vartheta_{n,\psi}^{\pmb{\lambda}}(\pmb{x})$,}\\
                 < 0 & \text{if $t>\vartheta_{n,\psi}^{\pmb{\lambda}}(\pmb{x})$},
             \end{cases}
           \end{align*}
           that is, for all $\pmb{x}\in X^n$, the value $\vartheta_{n,\psi}^{\pmb{\lambda}}(\pmb{x})$ is
           a point of sign change for the function $\psi_{\pmb{x},\pmb{\lambda}}$.
           If there is no confusion, instead of $\vartheta_{n,\psi}^{\pmb{\lambda}}$ we simply write $\vartheta_n^{\pmb{\lambda}}$.
          We may call $\vartheta_{n,\psi}^{\pmb{\lambda}}(\pmb{x})$
          as a weighted generalized $\psi$-estimator for some unknown parameter in $\Theta$ based
          on the realization $\bx=(x_1,\ldots,x_n)\in X^n$ and weights $(\lambda_1,\ldots,\lambda_n)\in\Lambda_n$.
   \end{enumerate}
\end{Def}

Given properties $[P_1], \ldots, [P_q]$ introduced in Definition~\ref{Def_Tn} (where $q\in\NN$), the subclass of $\Psi(X,\Theta)$ consisting of elements possessing the properties $[P_1],\ldots,[P_q]$,  will be denoted by $\Psi[P_1,\ldots,P_q](X,\Theta)$.

We call the attention to the fact that, given $\psi\in\Psi(X,\Theta)$ and $n\in\NN$, by Remark \ref{Rem_2}, the function $\Theta\ni t\mapsto \psi_{\pmb{x}}(t)$ can have at most one point of sign change for each $\bx\in X^n$.
Consequently, if $\psi\in\Psi[T_n](X,\Theta)$, then the generalized $\psi$-estimator $\vartheta_{n,\psi}(\pmb{x})$ introduced in part (ii) of Definition \ref{Def_Tn} is unique for each $\bx\in X^n$.
A similar conclusion holds for the weighted generalized $\psi$-estimator introduced in part (iv) of Definition \ref{Def_Tn}.
Therefore, in our forthcoming results (e.g., part (vi) of Theorem \ref{Thm_M_est_uniq}), when we establish that under some appropriate conditions a function $\psi\in\Psi(X,\Theta)$ satisfies
the property $[T_n]$ for each $n\in\NN$ or the property $[T_n^\blambda]$ for each $n\in\NN$ and $\blambda\in\Lambda_n$, then it means that our result in question provides conditions under which
the (weighted) generalized $\psi$-estimator exists uniquely.
If $\psi$ is continuous in its second variable as well, then such results provide conditions under which the usual $\psi$-estimator exists uniquely.

In the next proposition, we study the measurability of (weighted) generalized $\psi$-estimators, provided that $X$ is a measurable space.

\begin{Pro}\label{Pro_measurability}
Let $(X,\cX)$ be a measurable space, let $\Theta$ be a nonempty open interval of $\RR$, let $n\in\NN$, $\psi\in\Psi[Z_n](X,\Theta)$, and suppose that $\psi$ is measurable in its first variable. Then $\vartheta_{n,\psi}:X^n\to \Theta$ is measurable with respect to the sigma-algebras $\cX^n$ and $\cB(\Theta)$.
\end{Pro}

\begin{proof}
For all $r\in\Theta$, we have that
 \begin{align}\label{help_measurability}
  \begin{split}
   \vartheta_{n,\psi}^{-1}((-\infty,r))
      &= \big\{ (x_1,\ldots,x_n)\in X^n : \vartheta_{n,\psi}(x_1,\ldots,x_n)<r \big\}\\
      &= \Big\{ (x_1,\ldots,x_n)\in X^n : \sum_{i=1}^n \psi(x_i,r) < 0 \Big\},
   \end{split}
 \end{align}
where the second equality is a consequence of the property $[Z_n]$ of $\psi$. Further, for all $r\in\Theta$, the measurability of the mapping $X\ni x\mapsto \psi(x,r)$ implies the measurability of the mapping
 $X^n\ni (x_1,\ldots,x_n)\mapsto (\psi(x_1,r),\ldots,\psi(x_n,r))$ with respect to the sigma-algebras $\cX^n$ and $\cB(\Theta^n)$, and hence, using that the summation $\Theta^n\ni (t_1,\ldots,t_n)\mapsto t_1+\cdots+t_n$ is continuous,
 we have that $X^n\ni (x_1,\ldots,x_n)\mapsto \sum_{i=1}^n \psi(x_i,r)$ is measurable as well.
By \eqref{help_measurability}, it implies that $\vartheta_{n,\psi}^{-1}((-\infty,r))\in\cX$ for all $r\in\Theta$.
Since the sigma-algebra generated by the family $\{ (-\infty,r)\cap \Theta, r\in\Theta\}$ of intervals coincides with the Borel sigma-algebra $\cB(\Theta)$, we get the desired measurability of the generalized $\psi$-estimator $\vartheta_{n,\psi}$.
\end{proof}

As a consequence of Proposition \ref{Pro_measurability}, if $\xi_1,\ldots,\xi_n$ are random variables
 on a probability space $(\Omega,\cA,\PP)$, then
 $\vartheta_{n,\psi}(\xi_1,\ldots,\xi_n)$ is a random variable (measurable with respect to the sigma-algebras $\cA$ and $\cB(\Theta)$),
 i.e., it is a statistic in the language of mathematical statistics.
A similar statement to Proposition \ref{Pro_measurability} could be formulated for weighted generalized $\psi$-estimators as well.

Next, we present some basic facts about the properties $[T_n]$ and $[T_n^{\pmb{\lambda}}]$ given in Definition \ref{Def_Tn}, which can be easily checked.

\begin{Rem}\label{Rem_1}
(i) If $n\in\NN$ and $\psi\in\Psi[T_n](X,\Theta)$, then for each $x_1,\ldots,x_n\in X$,
 the equation \eqref{psi_est_equation} can have at most one solution.

(ii) If $n\in\NN$, $\psi\in\Psi[T_n](X,\Theta)$ and $\psi$ is continuous in its second variable, then $t=\vartheta_n(x_1,\dots,x_n)$ is the unique solution to \eqref{psi_est_equation}, and is called the $\psi$-estimator (Z-estimator) based on the observations $x_1,\ldots,x_n\in X$.
In particular, if $\psi\in\Psi[T_1](X,\Theta)$  and $\psi$ is continuous in its second variable, then, for each $x\in X$, the equation $\psi(x,t)=0$, $t\in\Theta$, has a unique solution $\vartheta_1(x)$.

(iii) If $\lambda_1=\cdots=\lambda_n>0$ with some $n\in\NN$, then $\vartheta_{n,\psi}^{(\lambda_1,\ldots,\lambda_n)} = \vartheta_{n,\psi}$.
\proofend
\end{Rem}

In the next remark, we point out an invariance property of the properties $[T_n]$ and $[T_n^{\pmb{\lambda}}]$
 given in Definition \ref{Def_Tn}.

\begin{Rem}
We introduce a notion of equivalence in $\Psi(X,\Theta)$ as follows.
We say that the maps $\psi,\varphi\in\Psi(X,\Theta)$ are \emph{equivalent} (denoted as $\psi\sim\varphi$)
 if there exists a positive function $h:\Theta\to\RR$ such that $\psi(x,t)=h(t)\varphi(x,t)$ is valid for all $(x,t)\in X\times \Theta$.
It is easy to see that $\sim$ is an equivalence relation on $\Psi(X,\Theta)$, furthermore,
 the properties $[T_n]$ and $[T_n^{\pmb{\lambda}}]$ are invariant with respect to this equivalence,
 that is, if $\psi\sim\varphi$ and $\varphi$ possesses the property $[T_n]$ (or the property $[T_n^{\pmb{\lambda}}]$),
 then $\psi$ also enjoys this property and $\vartheta_{n,\varphi}=\vartheta_{n,\phi}$
 (resp.\ $\vartheta_{n,\varphi}^{\pmb{\lambda}}=\vartheta_{n,\phi}^{\pmb{\lambda}}$).
\proofend
\end{Rem}

\begin{Def}\label{Def_level_of_increase}
Let $\Theta$ be a nonempty open interval of $\RR$ and $f:\Theta\to\RR$ be a function.
We say that $y\in\RR$ is a \emph{level of increase for $f$} if
$u,v\in\Theta$ and $f(v)\leq y\leq f(u)$ imply $v\leq u$.
\end{Def}

\begin{Rem}\label{Rem_level_of_increase}
(i) If $y\in\RR$ is such that either $f(u)>y$ for all $u\in\Theta$ or $f(u)<y$ for all $u\in\Theta$,
 then $y$ is automatically a level of increase for $f$.

(ii) If $y\in\RR$ is a level of increase for $f$, then the inverse image $f^{-1}(\{y\})$ is either empty or a singleton.
In general, the converse of the previous statement is not true.
To give a counterexample, let us consider the function $f:\Theta\to\RR$ given by $f(\vartheta):=1$ if $\vartheta\in\Theta$ is rational, and $f(\vartheta):=0$ if $\vartheta\in\Theta$ is irrational.
Then $f^{-1}(\{\frac{1}{2}\})=\emptyset$, but $\frac{1}{2}$ is not a level of increase for $f$.

(iii) Under the condition of Definition \ref{Def_level_of_increase}, $y\in\RR$ is a level of increase for $f$ if and only if the relations $u\in\Theta$ and $y\leq f(u)$ imply that $y<f(v)$ for all $v\in\Theta$ with $u<v$.
Indeed, if $y\in\RR$ is a level of increase for $f$, and $u,v\in\Theta$ are such that $u<v$ and $y\leq f(u)$, then
 $y<f(v)$, since otherwise $f(v)\leq y \leq f(u)$ would yield that $v\leq u$, leading us to a contradiction.
Conversely, assume that $y\in\RR$ is such that the relations $u\in\Theta$ and $y\leq f(u)$ imply that $y<f(v)$ for all $v\in\Theta$ with $u<v$.
If $u,v\in\Theta$ are such that $f(v)\leq y\leq f(u)$, then $v>u$ cannot hold, since it would yield that
 $y<f(v)$, leading us to a contradiction.
\proofend
\end{Rem}

In the following lemma, we establish a connection between the notions of point of sign change and level of increase.

\begin{Lem}\label{Lem_level_of_increase_and-sign_change}
Let $\Theta$ be a nonempty open interval of $\RR$,
 $f:\Theta\to\RR$ be a function, and $y\in\RR$.
Then $y$ is a level of increase for $f$ if and only if one of the following assertions holds:
\begin{enumerate}[(i)]\itemsep=-3pt
 \item $y<f$ on $\Theta$.
 \item $y>f$ on $\Theta$.
 \item There exists a point of sign change for the function $y-f$.
\end{enumerate}
\end{Lem}

The proof of Lemma \ref{Lem_level_of_increase_and-sign_change} and that of all the forthcoming results in this section can be found in Section \ref{Section_Proof_for_Section2}.

\begin{Lem}\label{Lem_level_of_increase}
Let $\Theta$ be a nonempty open interval and $f:\Theta\to\RR$ be a function.
If the levels of increase for $f$ form a dense subset in the convex hull of $f(\Theta)$, then $f$ is increasing.
The function $f$ is strictly increasing if and only if every element of $f(\Theta)$ is a level of increase for $f$.
Furthermore, if $g:H\to\RR$ is a strictly increasing function, where $H$ is a set containing $f(\Theta)$,
 and $y\in H$ is level of increase for $f$, then $g(y)$ is a level of increase for $g\circ f$.
\end{Lem}

We recall also a definition due to P\'ales \cite{Pal2003}: given a nonempty open interval $\Theta$ of $\RR$ and $\vare>0$, a function $f:\Theta\to \RR$ satisfying the inequality $f(u)\leq f(v) +\vare$ for all $u<v$, $u,v\in \Theta$ is called \emph{$\vare$-increasing}.
If the inequality is strict for all $u<v$, $u,v\in \Theta$, then $f$ is said to be \emph{strictly $\vare$-increasing}.
We note that P\'ales \cite[Theorem 3]{Pal2003} offers the following simple characterization:
 a function $f:\Theta\to\RR$ is $\vare$-increasing if and only if there exists an increasing function $g:\Theta\to\RR$
 such that $\|f-g\|_\infty:=\sup_{u\in \Theta}| f(u) - g(u)|\leq \frac{\vare}{2}$.

The next lemma describes a connection between levels of increase for a function $f:\Theta\to \RR$ and its $\vare$-increasingness property.

\begin{Lem}\label{Lem_e_monoton}
Let $\Theta$ be a nonempty open interval, let $n\in\NN$ and let $y_0<\dots<y_n$ be real numbers.
Assume that $y_0,\dots,y_{n-1}$ are levels of increase for a function $f:\Theta\to\RR$ and $f(\Theta)\subseteq[y_0,y_n]$.
Then $f$ is strictly $\varepsilon$-increasing with $\varepsilon:=\max\{y_1-y_0,\dots,y_n-y_{n-1}\}$.
\end{Lem}

Now, we state our first main result by presenting necessary as well as sufficient conditions for the properties $[T_n]$ and $[T_n^{\pmb{\lambda}}]$.

\begin{Thm}\label{Thm_M_est_uniq}
Let $X$ be a nonempty set, $\Theta$ be a nonempty open interval of $\RR$, and $ \psi\in\Psi[T_1](X,\Theta)$.
\begin{enumerate}[(i)]\itemsep=-2pt
 \item If $\psi\in\Psi[T_2^{(\lambda_1,\lambda_2)}](X,\Theta)$ for some $(\lambda_1,\lambda_2)\in(0,\infty)^2$, then, for each $x,y\in X$ with $\vartheta_1(x)<\vartheta_1(y)$, the numbers $\frac{\lambda_1}{\lambda_2}$ and $\frac{\lambda_2}{\lambda_1}$ are levels of increase for the function
 \begin{align}\label{function_newhanyados0}
   (\vartheta_1(x),\vartheta_1(y))\ni t \mapsto  -\frac{\psi(x,t)}{\psi(y,t)}.
 \end{align}
 \item If $\psi\in\Psi[T_n^{(\lambda_1,\ldots,\lambda_n)}](X,\Theta)$ for some $n\in\NN\setminus\{1\}$ and $(\lambda_1,\ldots,\lambda_n)\in(0,\infty)^n$, then, for each $x,y\in X$ with $\vartheta_1(x)<\vartheta_1(y)$, the numbers $\frac{\lambda_1+\cdots+\lambda_k}{\lambda_{k+1}+\cdots+\lambda_n}$ and
 $\frac{\lambda_{k+1}+\cdots+\lambda_n}{\lambda_1+\cdots+\lambda_k}$, $k\in\{1,\ldots,n-1\}$, are levels of increase for the function \eqref{function_newhanyados0}.
 \item If $\psi\in\Psi[T_n](X,\Theta)$ for some $n\in\NN\setminus\{1\}$, then, for each $x,y\in X$ with $\vartheta_1(x)<\vartheta_1(y)$, the elements of the set $\{\frac{k}{n-k}\mid k\in\{1,\dots,n-1\}\}$ are levels of increase for the function \eqref{function_newhanyados0}.
 \item If $\psi\in\Psi[T_n](X,\Theta)$ for infinitely many $n\in\NN$, then for each $x,y\in X$ with $\vartheta_1(x)<\vartheta_1(y)$, the function \eqref{function_newhanyados0} is increasing. In addition, if for each $m\in\NN$ there exists $n\in\NN$ such that $m$ divides $n$ and $\psi\in\Psi[T_n](X,\Theta)$, then, for each $x,y\in X$ with $\vartheta_1(x)<\vartheta_1(y)$, every positive rational number is a level of increase for the function \eqref{function_newhanyados0}.

 \item If $\psi\in\Psi[T_2^{\pmb{\lambda}}](X,\Theta)$ for each $\pmb{\lambda}\in\Lambda_2$, then for each $x,y\in X$ with $\vartheta_1(x)<\vartheta_1(y)$, the function \eqref{function_newhanyados0} is strictly increasing.
\item If $\psi\in\Psi[Z_1](X,\Theta)$ and, for each $x,y\in X$ with $\vartheta_1(x)<\vartheta_1(y)$, the function \eqref{function_newhanyados0} is strictly increasing, then $ \psi\in\Psi [T_n^{\pmb{\lambda}}](X,\Theta)$ for each $n\in\NN$ and $\pmb{\lambda}\in\Lambda_n$.
\end{enumerate}
\end{Thm}

We emphasize that, for $\psi\in\Psi[T_1](X,\Theta)$, the assertion (vi) of Theorem \ref{Thm_M_est_uniq} provides a sufficient condition for the existence and uniqueness of a weighted generalized $\psi$-estimator.
If, in addition, $\psi$ is continuous in its second variable, then it gives a sufficient condition for the existence and uniqueness of a (usual) $\psi$-estimator.

The following statement establishes three equivalent conditions under the property $[Z_1]$.

\begin{Cor}\label{Cor_M_est_uniq}
Let $X$ be a nonempty set, $\Theta$ be a nonempty open interval of $\RR$, and $\psi\in\Psi[Z_1](X,\Theta)$.
Then the following assertions are equivalent:
\begin{enumerate}[(i)]\itemsep=-2pt
\item For each $x,y\in X$ with $\vartheta_1(x)<\vartheta_1(y)$, the function \eqref{function_newhanyados0} is strictly increasing.

\item For each $\pmb{\lambda}\in\Lambda_2$, we have $\psi\in\Psi[T_2^{\pmb{\lambda}}](X,\Theta)$.
\item For each $n\in\NN$ and $\pmb{\lambda}\in\Lambda_n$, we have $\psi\in\Psi[T_n^{\pmb{\lambda}}](X,\Theta)$.
\end{enumerate}
\end{Cor}

In part (ii) of the next proposition, we provide a sufficient condition (which does not involve the property $[Z_1]$)
 under which $\psi$ has the property $[T_n^{\pmb{\lambda}}]$ for each $n\in\NN$ and $\pmb{\lambda}\in\Lambda_n$.
We also call the attention to the fact that our proof is elementary in the sense that it does not depend on Theorem \ref{Thm_M_est_uniq}.

\begin{Pro}\label{Pro_M_est_uniq}
Let $X$ be a nonempty set, $\Theta$ be a nonempty open interval of $\RR$, and
$\psi\in\Psi[T_1](X,\Theta)$.
\begin{enumerate}[(i)]\itemsep=-2pt
 \item If for each $x\in X$, the function $\Theta\ni t\mapsto \psi(x,t)$ is (strictly) decreasing, then for each $x,y\in X$ with $\vartheta_1(x)<\vartheta_1(y)$, the
       function \eqref{function_newhanyados0}  is (strictly) increasing.
 \item If for each $x\in X$, the function $\Theta\ni t\mapsto \psi(x,t)$ is strictly decreasing, then $\psi\in\Psi[T_n^{\pmb{\lambda}}](X,\Theta)$ for each $n\in\NN$ and $\pmb{\lambda}\in\Lambda_n$.
\end{enumerate}
\end{Pro}

The next proposition establishes a connection between the property $[T_n]$ and the strictly $\frac1n$-increasingness
of some appropriately defined functions.

\begin{Pro}\label{Pro_3}
Let $X$ be a nonempty set, $\Theta$ be a nonempty open interval of $\RR$, and $\psi\in\Psi[T_1](X,\Theta)$.
If $\psi\in\Psi[T_n](X,\Theta)$ for some $n\in\NN\setminus\{1\}$, then, for each $x,y\in X$ with $\vartheta_1(x)<\vartheta_1(y)$, the function
 \begin{align}\label{function_newhanyados}
   (\vartheta_1(x),\vartheta_1(y))\ni t \mapsto  \frac{\psi(x,t)}{\psi(x,t)-\psi(y,t)}
 \end{align}
 is strictly $\frac1n$-increasing.
\end{Pro}

The following two results describe the hierarchy among the properties $([T_n])_{n\in\NN}$ and establish a kind of 'grouping' property of
 the property $[T_n^{\pmb{\lambda}}]$.

\begin{Pro}\label{T_n}
Let $X$ be a nonempty set, $\Theta$ be a nonempty open interval of $\RR$, and $\psi\in\Psi(X,\Theta)$.
If $\psi\in\Psi[T_n](X,\Theta)$ for some $n\in\NN$, then $\psi\in\Psi[T_m](X,\Theta)$ for any $m\in\{1,\dots,n\}$ that divides $n$.
\end{Pro}

\begin{Pro}\label{Pro_T_n_lambda}
Let $X$ be a nonempty set, $\Theta$ be a nonempty open interval of $\RR$, and
$\psi\in\Psi[T_n^{\pmb{\lambda}}](X,\Theta)$ for some $n\in\NN$ and $\pmb{\lambda}=(\lambda_1,\ldots,\lambda_n)\in\Lambda_n$.
Let $m\in\{1,\ldots,n\}$ and $H_1,\dots,H_m$ be nonempty pairwise disjoint subsets of $\{1,\ldots,n\}$ such that $H_1\cup\cdots \cup H_m=\{1,\ldots,n\}$.
For each $\alpha\in \{1,\ldots, m\}$, define $\mu_\alpha:=\sum_{i\in H_\alpha}\lambda_i$. Then
$\pmb{\mu}:=(\mu_1,\ldots,\mu_m)\in \Lambda_m$ and $\psi\in\Psi[T_m^{\pmb{\mu}}](X,\Theta)$.
\end{Pro}

Note that Proposition \ref{Pro_T_n_lambda} implies Proposition \ref{T_n}.
Indeed, let $\psi\in\Psi[T_n](X,\Theta)$ for some $n\in\NN$.
If $m$ is a divisor of $n$, then $n=km$ with some $k\in\NN$, and hence $H_j:=\{(j-1)k+1,\ldots,jk\}$, $j=1,\ldots,m$, are nonempty pairwise disjoint subsets such that $H_1\cup\cdots \cup H_m=\{1,\ldots,n\}$.
With the choice $\blambda:=(\lambda_1,\ldots,\lambda_n):=(1,\ldots,1)\in \Lambda_n$,
we have  $\psi$ is a $T_n^{\pmb{\lambda}}$-function and $\mu_\alpha=k$, $\alpha\in\{1,\ldots,m\}$.
Hence Proposition \ref{Pro_T_n_lambda} yields that $\psi$ is a $T_m^{\pmb{\mu}}$-function.
Since $\mu_\alpha$, $\alpha\in\{1,\ldots, m\}$, are all the same positive constant $k$, we have
 $\psi$ is a $T_m$-function as well.
Further, we note that Proposition \ref{Pro_T_n_lambda} may be useful to see, for example, that the property $[T_3^{(\lambda_1,\lambda_2,\lambda_3)}]$ of $\psi$ implies the property $[T_2^{(\lambda_1+\lambda_2,\lambda_3)}]$ of $\psi$, where $(\lambda_1,\lambda_2,\lambda_3)\in\Lambda_3$.

The next example demonstrates that the property $[T_2]$ can already fail to hold for a $T_1$-function.

\begin{Ex}\label{Ex_T2_fail}
Let $m\in\NN$, $X:=\{x_1,\dots,x_m\}$, $\Theta:=\RR$ and let $w_1,\dots,w_m>0$. Define $\psi:X\times\Theta\to\RR$ by
\[
 \psi(x_i,t):=\begin{cases}
                 w_i & \text{if $t<i$,}\\
                 -w_i & \text{if $t\geq i$,}
               \end{cases}
               \qquad i\in\{1,\ldots,m\}.
\]
Then $\psi\in\Psi[T_1](X,\Theta)$ and $\vartheta_1(x_i)=i$ holds for all $i\in\{1,\dots,m\}$.
One can easily see that the property $[T_2]$ holds if and only if we have $w_i\neq w_j$ for all distinct $i,j\in\{1,\dots,m\}$,
Indeed, if $1\leq i<j\leq m$, then
\[
\psi(x_i,t) + \psi(x_j,t)
       =\begin{cases}
          w_i+w_j>0 & \text{if $t<i$,}\\
          -w_i+w_j & \text{if $i\leq t<j$,}\\
          -w_i-w_j<0 & \text{if $j\leq t$.}
        \end{cases}
\]
This function has a point of sign change if and only if $w_i\neq w_j$, as desired.
We also get that if $\psi\in\Psi[T_2](X,\Theta)$, then
 \[
 \vartheta_2(x_i,x_j)=\vartheta_2(x_j,x_i)
 =\begin{cases}
    i & \text{if $w_i>w_j$,} \\
    j & \text{if $w_i<w_j$.}
  \end{cases}
 \]
Furthermore, $\vartheta_2(x_i,x_i)=\vartheta_1(x_i)=i$ for all $i\in\{1,\dots,m\}$.
\proofend
\end{Ex}

In what follows, we give an example to point out that in part (vi) of Theorem \ref{Thm_M_est_uniq}
 the assumption that $\psi(x,\vartheta_1(x))=0$, $x\in X$, cannot be omitted.

\begin{Ex}\label{Ex_not_omitted}
Let $X:=\{x_1,x_2\}$ (with $x_1\neq x_2$) and $\Theta:=\RR$.
Let
 \[
  \psi(x_1,t):=\begin{cases}
                 2 & \text{if $t<1$,}\\
                 -t & \text{if $1\leq t\leq 2$,}\\
                 -2 & \text{if $t>2$,}\\
               \end{cases}
 \]
 and
 \[
  \psi(x_2,t):=\begin{cases}
                 1 & \text{if $t<2$,}\\
                 2 & \text{if $t=2$,}\\
                 -1 & \text{if $t>2$.}\\
               \end{cases}
 \]
Then $\psi\in\Psi[T_1](X,\Theta)$ with $\vartheta_1(x_1)=1$ and $\vartheta_1(x_2)=2$, and $\psi(x_i,\vartheta_1(x_i))\neq0$ for $i\in\{1,2\}$.
We also note that the function
 \[
   (\vartheta_1(x_1),\vartheta_1(x_2))=(1,2)\ni t\mapsto \frac{\psi(x_1,t)}{\psi(x_2,t)} = -t
 \]
 is strictly decreasing.
 However, $\psi$ is not a $T_2$-function, since
 \[
   \psi(x_1,t) + \psi(x_2,t)
       =\begin{cases}
          2+1=3>0 & \text{if $t<1$,}\\
          -1+1=0 & \text{if $t=1$,}\\
          -t+1<0 & \text{if $1<t<2$,}\\
          -2+2=0 & \text{if $t=2$,}\\
          -2-1=-3<0 & \text{if $t>2$},
        \end{cases}
 \]
which shows that $\RR\ni t\mapsto \psi(x_1,t) + \psi(x_2,t)$ does not have a point of sign change.
\proofend
\end{Ex}

In what follows, as an application of Proposition \ref{Pro_M_est_uniq},
 we present an example of a large class of functions $\psi:X\times \Theta\to\RR$,
 which possesses the property $[T_n^{\blambda}]$ for each $n\in\NN$ and $\blambda\in\Lambda_n$ and
 for which the point of sign change $\vartheta_{n,\psi}^{\blambda}(\bx)$ (where $\bx\in X^n$)
 has an explicit form.
This class of functions may be called the class of Bajraktarevi\'c-type functions, motivated by the representation
 of Bajraktarevi\'c means as special deviation means.
For the description of this class of functions, we need to recall the notion of generalized left inverse of a strictly monotone (but not necessarily continuous) function defined on a nonempty open
 interval of $\RR$, see, e.g., Gasi\'nski and Papageorgiou \cite[Proposition 1.55 and the subsequent comment]{GasPap}
  and  Gr\"unwald and P\'ales \cite[Lemma 1]{GruPal}.
The notion of generalized left inverse in question is likely to be well-known and its properties are established for a while,
 but we could not trace the roots, and therefore we refer to the recent treatments appearing in \cite{GasPap} and \cite{GruPal}.

\begin{Lem}\label{Lem_gen_left_inverse}
Let $\Theta$ be a nonempty open interval of $\RR$, let $f:\Theta\to\RR$ be a strictly increasing function.
Then there exists a uniquely determined monotone function $g:\conv(f(\Theta))\to \Theta$ such that $g$ is the left inverse of $f$, i.e.,
\[
   (g\circ f)(x)=x, \qquad x\in \Theta.
\]
Furthermore, $g$ is monotone in the same sense as $f$, is continuous, and the following relation holds:
\[
  (f\circ g)(y)=y, \qquad y\in f(\Theta).
\]
\end{Lem}

The function $g:\conv(f(\Theta))\to \Theta$ described in Lemma \ref{Lem_gen_left_inverse} is called the \emph{generalized left inverse of the strictly increasing function $f:\Theta\to\RR$} and is denoted by $f^{(-1)}$.
In fact, by the proof of Lemma 1 in Gr\"unwald and P\'ales \cite{GruPal}, it also turns out that
 \[
   g(y) = \sup\{ u\in\Theta : f(u)\leq y\} = \inf\{ u\in\Theta : f(u)\geq y\},\qquad y\in\conv(f(\Theta)).
 \]
It is clear that the restriction of $f^{(-1)}$ to $f(\Theta)$ is the inverse of $f$ in the standard sense.
Therefore, $f^{(-1)}$ is the continuous and monotone extension of the inverse of $f$ to the smallest interval
 containing the range of $f$, that is, to the convex hull of $f(\Theta)$.

\begin{Def}\label{Ex_Bajrak_emp}
Let $X$ be a nonempty set, $\Theta$ be a nonempty open interval of $\RR$, $f:\Theta\to\RR$ be a strictly increasing function, and $p:X\to\RR_{++}$ and $\varphi:X\to\conv(f(\Theta))$ be functions.
In terms of these functions, define $\psi:X\times\Theta\to\RR$ by
 \begin{align}\label{help16}
   \psi(x,t):=p(x)(\varphi(x)-f(t)), \qquad x\in X, \; t\in\Theta.
 \end{align}
The function $\psi$ defined by \eqref{help16} is said to be of Bajraktarevi\'c-type. In particular, if $p=1$ is a constant function, then $\psi$ is said to be of quasi-arithmetic-type.
\end{Def}

\begin{Pro}\label{Pro_BajType_est}
Under the assumptions of Definition~\ref{Ex_Bajrak_emp}, we have that $\psi\in\Psi[T_n^{\blambda}](X,\Theta)$ for each $n\in\NN$ and $\blambda=(\lambda_1,\dots,\lambda_n)\in\Lambda_n$, and
 \begin{align}\label{help14}
  \vartheta_{n,\psi}^{\blambda}(\bx)
  =f^{(-1)}\bigg(\frac{\lambda_1p(x_1)\varphi(x_1)+\dots+\lambda_np(x_n)\varphi(x_n)}{\lambda_1p(x_1)+\dots+\lambda_np(x_n)}\bigg)
 \end{align}
 for all $\bx=(x_1,\dots,x_n)\in X^n$.
In particular, the equality $\vartheta_{1,\psi}=f^{(-1)}\circ\varphi$ holds.
\end{Pro}

One may call the value $\vartheta_{n,\psi}^{\blambda}(\bx)$ given by \eqref{help14} as a Bajraktarevi\'c-type $\psi$-estimator of some unknown parameter in $\Theta$ based on the realization $\bx=(x_1,\ldots,x_n)\in X^n$ and weights $\blambda=(\lambda_1,\ldots,\lambda_n)\in\Lambda_n$ corresponding to the Bajraktarevi\'c-type function given by \eqref{help16}.
In particular, if $p=1$ is a constant function in \eqref{help16}, then we speak about a quasi-arithmetic-type $\psi$-estimator.

Note that in case of $X:=\Theta$ and $\varphi:=f$,  Proposition \ref{Pro_BajType_est} reduces to Theorem 3 in Gr\"unwald and P\'ales \cite{GruPal} for Bajraktarevi\'c means.
In addition, if $p=1$ is a constant function, then Proposition~\ref{Pro_BajType_est} is about generalized quasi-arithmetic means (here, we use the term 'generalized', since for usual quasi-arithmetic means, the function $f$ is not only strictly increasing, but continuous as well).

In the next example, we point out that, in general, one cannot omit the restriction that $m$ divides $n$ in Proposition \ref{T_n}.
For another example, see the case of empirical median discussed in Example \ref{Ex_median}.

\begin{Ex}\label{Ex_div_noomit}
Let $X$ be an arbitrary set with at least two  distinct elements and let $X_1,X_2\subseteq X$ be nonempty disjoint subsets such that $X_1\cup X_2=X$ and let $w_1,w_2>0$.
Define $\psi:X\times\Theta\to\RR$ by
\[
 \psi(x,t):=\begin{cases}
                 w_i & \text{$t<i$,}\\
                 -w_i & \text{$t\geq i$,}
               \end{cases} \qquad \mbox{if $x\in X_i$}.
\]
Then $\psi\in\Psi[T_1](X,\Theta)$ and, for $i\in\{1,2\}$ and $x\in X_i$, we have $\vartheta_1(x)=i$.

Let $n,k\in\NN$ such that $k$ is not a divisor of $n$ (which implies that $k\geq 2$). Then
\begin{equation}\label{eq:1/k}
    \frac{1}{k}
   \not\in\bigg\{\frac{1}{n},\dots,\frac{n-1}{n}\bigg\}.
\end{equation}
Indeed, on the contrary, if the inclusion were valid, then $1/k$ would be of the form $m/n$ for some $m\in\{1,\dots,n-1\}$, yielding that $n=km$, which contradicts the assumption $k\nmid n$.

Assuming that $w_1=k-1$ and $w_2=1$, we prove that $\psi$ is a $T_n$-function, but it is not a $T_k$-function.

To show that $\psi$ is a $T_n$-function, let $\pmb{y}:=(y_1,\dots,y_n)\in X^n$.
For $j\in\{1,2\}$, define the set $S_j:=\{i\in\{1,\dots,n\}:y_i\in X_j\}$.
Then $\{S_1,S_2\}$ forms a partition of $\{1,\dots,n\}$.
Let $n_j$ denote the cardinality of $S_j$, $j\in\{1,2\}$.
Then $n=n_1+n_2$ and
 \[
 \psi_{\pmb{y}}(t):=\sum_{i=1}^n\psi(y_i,t)
 =\begin{cases}
          n_1w_1+n_2w_2>0 & \text{if $t<1$,}\\
          -n_1w_1+n_2w_2 & \text{if $1\leq t<2$,}\\
          -n_1w_1-n_2w_2<0 & \text{if $2\leq t$.}
    \end{cases}
\]
Using condition \eqref{eq:1/k} and $k\geq2$, we have
\[
-n_1w_1+n_2w_2=-n_1(k-1)+(n-n_1)=n-n_1k\neq0.
\]
Therefore, the point of sign change for the function $\psi_{\pmb{y}}$ equals $1$ if $-n_1w_1+n_2w_2<0$ and equals $2$ if $-n_1w_1+n_2w_2>0$.
This proves that $\psi\in\Psi[T_n](X,\Theta)$.

To verify that $\psi$ is not a $T_k$-function, let $x_1\in X_1$ and $x_2\in X_2$ be fixed and let $\pmb{z}:=(x_1,x_2,\dots,x_2)\in X^k$. Then
\[
 \psi_{\pmb{z}}(t):=
  \psi(x_1,t)+(k-1)\psi(x_2,t)
 =\begin{cases}
          w_1+(k-1)w_2>0 & \text{if $t<1$,}\\
          -w_1+(k-1)w_2=0 & \text{if $1\leq t<2$,}\\
          -w_1-(k-1)w_2<0 & \text{if $2\leq t$.}
    \end{cases}
\]
Therefore, the function $\psi_{\pmb{z}}$ does not have a point of sign change, and, consequently, $\psi$ is not a $T_k$-function, as desired.
\proofend
\end{Ex}

\section{Existence and uniqueness of the point of sign change of $\psi$-expectation functions}
\label{Section_Ex_uniq_theoretical}

In this section, we investigate Problem 2 presented in the Introduction.

As it was mentioned in the Introduction, in case of simple random variables Problem 2 is a special case of Problem 1.
More precisely, if $\psi\in\Psi(X,\Theta)$ and $\xi$ is a simple random variable such that $\PP(\xi=x_i)=p_i$, $i=1,\ldots,n$, where $n\in\NN$, $(x_1,\ldots,x_n)\in X^n$ and $p_1,\ldots,p_n\geq 0$, $p_1+\cdots+p_n=1$,
 then $\EE(\psi(\xi,t)) = \sum_{i=1}^n p_i\psi(x_i,t)$, $t\in\Theta$.
In addition, if $\psi$ is a $T_n^{(p_1,\ldots,p_n)}$-function, then, by definition, the function \eqref{fgv_elmeleti}
 has a unique point of sign change.
Further, Theorem \ref{Thm_M_est_uniq} provides some necessary as well as some sufficient conditions
under which $\psi$ possesses the property $[T_n^{(p_1,\ldots,p_n)}]$.
In case of a general (not necessarily discrete) random variable $\xi$ and $\psi\in\Psi(X,\Theta)$,
in our forthcoming results Theorem \ref{Thm_theor_M_est_uniq} and Proposition \ref{Pro_theor_M_est_uniq2}, we derive sufficient conditions on $\xi$ and $\psi$ under which there exists a unique point of sign change of the corresponding $\psi$-expectation function.

Next, we present our second main result in which we give a set of sufficient conditions in order that the function given by \eqref{fgv_elmeleti} have a unique point of sign change.

\begin{Thm}\label{Thm_theor_M_est_uniq}
Let $(X,\cX)$ be a measurable space, $\Theta$ be a nonempty open interval of $\RR$, $\psi:X\times \Theta\to\RR$ be a function,
 and $\xi:\Omega\to X$ be a random variable defined on a probability space $(\Omega,\cA,\PP)$.
Let us suppose that
 \begin{enumerate}[(i)]
  \item $\psi\in\Psi[Z_1](X,\Theta)$,
  \item for each $x,y\in X$ with $\vartheta_1(x)<\vartheta_1(y)$, the function \eqref{function_newhanyados0} is strictly increasing,
  \item $\psi$ is measurable in its first variable,
  \item $\EE(\vert\psi(\xi,t)\vert)<\infty$ for each $t\in\Theta$,
  \item there exist $s_0,t_0\in\Theta$ such that
  $\EE(\psi(\xi,s_0))\geq0$ and $\EE(\psi(\xi,t_0))\leq0$.
 \end{enumerate}
Then the map $\Theta\ni t\to\EE(\psi(\xi,t))$ admits a unique point of sign change in $\Theta$.
\end{Thm}

The proof of Theorem \ref{Thm_theor_M_est_uniq} and that of all the forthcoming results
in this section can be found in Section \ref{Section_Proof_for_Section3}.

Next, we provide a set of sufficient conditions (which does not involve the condition $\psi(x,\vartheta_1(x))=0$ for each $x\in X$) under which the map $\Theta\ni t\to\EE(\psi(\xi,t))$ also has a unique point of sign change.

\begin{Pro}\label{Pro_theor_M_est_uniq2}
Let $(X,\cX)$ be a measurable space, $\Theta$ be a nonempty open interval of $\RR$, $\psi:X\times \Theta\to\RR$ be a function,
 and $\xi:\Omega\to X$ be a random variable defined on a probability space $(\Omega,\cA,\PP)$.
Let us suppose that
 \begin{enumerate}[(i)]
  \item for each $x\in X$, the function $\Theta\ni t\mapsto \psi(x,t)$ is strictly decreasing,
  \item $\psi$ is measurable in its first variable,
  \item $\EE(\vert\psi(\xi,t)\vert)<\infty$ for each $t\in\Theta$,
  \item there exist $s_0,t_0\in\Theta$ such that
  $\EE(\psi(\xi,s_0))\geq0$ and $\EE(\psi(\xi,t_0))\leq0$.
 \end{enumerate}
Then the function $\Theta\ni t\to\EE(\psi(\xi,t))$ admits a unique point of sign change in $\Theta$.
\end{Pro}

Next, we formulate a corollary of Theorem \ref{Thm_theor_M_est_uniq}, which is in fact part (vi) of Theorem \ref{Thm_M_est_uniq}.

\begin{Cor}\label{Cor_M_est_uniq+}
Let $X$ be a nonempty set, $\Theta$ be a nonempty open interval of $\RR$, and $\psi\in\Psi[Z_1](X,\Theta)$.
If, for each $x,y\in X$ with $\vartheta_1(x)<\vartheta_1(y)$, the function \eqref{function_newhanyados0} is strictly increasing, then $\psi\in\Psi[T_n^{\blambda}](X,\Theta)$ for each $n\in\NN$ and $\blambda\in\Lambda_n$.
\end{Cor}

In what follows, we present a particular case of Proposition \ref{Pro_theor_M_est_uniq2}, which can be considered as a counterpart of Proposition \ref{Pro_BajType_est} for Bajraktarevi\'c-type functions $\psi$.

\begin{Pro}\label{Pro_Bajrak_theor}
Let $(X,\cX)$ be a measurable space, $\Theta$ be a nonempty open interval of $\RR$, $f:\Theta\to\RR$ be a strictly increasing function, and $p:X\to\RR_{++}$ and $\varphi:X\to \conv(f(\Theta))$ be measurable functions. Define $\psi:X\times\Theta\to\RR$ by \eqref{help16}. Further, let  $(\Omega,\cA,\PP)$ be a probability space, $\xi:\Omega\to X$ be a random variable such that $\EE(p(\xi)|\varphi(\xi)|)<\infty$ and $\EE(p(\xi))<\infty$.%
Then the function $\Theta\ni t\to\EE(\psi(\xi,t))$  admits a unique point of sign change in $\Theta$ which is given by
 \[
   f^{(-1)}\bigg(\frac{\EE(p(\xi)\varphi(\xi))}{\EE(p(\xi))}\bigg).
 \]
\end{Pro}


The following auxiliary result is instrumental for the proof of Proposition \ref{Pro_Bajrak_theor}.

\begin{Lem}\label{Lem_conv_hull}
Let $(X,\cX)$ be a measurable space and  $p:X\to \RR_{++}$ and $\varphi:X\to \RR$ be measurable functions.
Further, let $\xi:\Omega\to X$ be a random variable on a probability space $(\Omega,\cA,\PP)$ such that $\EE(p(\xi))<\infty$ and $\EE(p(\xi)|\varphi(\xi)|)<\infty$.
Then
  \[
    \frac{\EE(p(\xi)\varphi(\xi))}{\EE(p(\xi))}\in \conv(\varphi(X)).
 \]
\end{Lem}

In the next example, we consider a particular form of $\psi$ which has been recently investigated by Mathieu \cite{Mat}: namely, let $\psi:\RR\times\RR\to\RR$,
 \begin{align}\label{help_psi_mathieu_new}
 \psi(x,t):=\sign(x-t)f(\vert x-t\vert),\qquad x,t\in\RR,
 \end{align}
where $f:\RR_+\to\RR_+$.
Mathieu \cite[Lemma 2]{Mat} has derived some sufficient conditions on $f$ and $\xi$ under which the equation
$\EE(\psi(\xi,t))=0$, $t\in\RR$, has a unique solution, for more details and our new results in this special case,
see the next example and Proposition \ref{Pro_Math_type}, respectively.

\begin{Ex}\label{Ex_Mat_elmeleti}
Let $X:=\RR$, $\Theta:=\RR$ and $\psi:\RR\times\RR\to\RR$ be given by \eqref{help_psi_mathieu_new}.
Given a random variable $\xi$, Mathieu \cite{Mat} has recently considered the problem of finding a
 unique element $t_0\in\Theta$ such that $\EE(\psi(\xi,t_0))=0$ holds, where $\psi$ has the form given in \eqref{help_psi_mathieu_new}
 such that $f$ admits the following properties (called Assumption 2 in Mathieu \cite{Mat}):
 \begin{itemize}
   \item[(a)] $f$ is continuous and differentiable Lebesgue almost everywhere,
   \item[(b)] $f(0)=0$,
   \item[(c)] $f$ is concave,
   \item[(d)] there exist $\beta,\gamma>0$ such that $\gamma \bone_{\{x\leq \beta\}}\leq f'(x)\leq 1$ Lebesgue a.e. $x\geq 0$.
 \end{itemize}
For historical fidelity, we note that Mathieu \cite{Mat} investigated a more general setup,
 he considered a random variable $\xi$ having values in a Hilbert space $\cH$,
 and a function $\psi:\cH\times\cH\to\cH$, $\psi(x,t):=\frac{x-t}{\Vert x-t \Vert} f(\Vert x-t \Vert)$ for $x\ne t$, $x,t\in\cH$,
 where $\Vert\cdot\Vert$ is the norm of the Hilbert space $\cH$ (the value of $\psi$ at $(x,x)$, $x\in\cH$, was not specified in Mathieu \cite{Mat}).

Mathieu \cite[Lemma 2]{Mat} has shown that (formulating his result only in the case of $\cH=\RR$) if $f$ admits the properties (a)-(d), $\EE(\vert \xi\vert)<\infty$
 and the inequality $\EE(\varrho(\vert \xi - \EE(\xi)\vert)) < \varrho(\beta)$ holds, where
 $\varrho(x):=\int_0^x f(u)\,\dd u$, $x\in\RR_+$, then there exists a unique element $t_0\in\Theta$ such that
 $\EE(\psi(\xi,t_0))=0$ holds.
Mathieu \cite{Mat} also noted that the assumptions under which existence and uniqueness of a solution in question was established are not the minimal ones, but he has not searched for possible minimal assumptions.
\proofend
\end{Ex}

Note that one can rewrite $\psi$ given by \eqref{help_psi_mathieu_new} as $\psi(x,t)=\widetilde{f}(x-t)$, $x,t\in\RR$,
 where $\widetilde{f}:\RR\to\RR$ denotes the odd extension of $f:\RR_+\to\RR_+$ to $\RR$, which is given by
 \begin{align*}
  \widetilde{f}(z)
   :=\begin{cases}
      f(z) & \mbox{if } z>0,\\
         0 & \mbox{if } z=0, \\
      -f(-z) & \mbox{if } z<0.
   \end{cases}
 \end{align*}

As a new result, we have the following proposition.


\begin{Pro}\label{Pro_Math_type}
If $f:\RR_+\to\RR_+$ is continuous and strictly increasing with $f(0)=0$ and $\lim_{z\to\infty} f(z)\in(0,\infty)$, then,
 for any random variable $\xi$, we have that $\EE(|\psi(\xi,t)|)<\infty$, $t\in\RR$, and the equation
 \[
   \EE(\psi(\xi,t))=\EE(\widetilde{f}(\xi-t))=0
 \]
has a unique solution with respect to $t\in\RR$.
\end{Pro}

Now, we compare the assumptions of Lemma 2 in Mathieu \cite{Mat} and those of Proposition~\ref{Pro_Math_type}.
Note that if a function $f:\RR_+\to\RR_+$ admits the properties (a)--(d) of Example \ref{Ex_Mat_elmeleti}, then it is not necessarily strictly increasing
 (for example, it may happen that $f(x)=f(\beta)$ for $x\geq \beta$, see, e.g., the Huber function \eqref{func_Huber}),
 so we cannot say that Proposition~\ref{Pro_Math_type} is a generalization of Lemma 2 in Mathieu \cite{Mat}.
However, the conditions of Proposition~\ref{Pro_Math_type} might be checked more easily than those of Lemma 2 in Mathieu \cite{Mat}
 in order to prove that the equation $\EE(\psi(\xi,t))=0$ have a unique solution with respect to $t\in\RR$.
For example, if $f:\RR_+\to\RR$, $f(z):=z/\sqrt{1+z^2/2}$, $z\in\RR_+$, then $f$ is a continuous and strictly increasing function
 starting from $0$ and $\lim_{z\to\infty}f(z)=\sqrt{2}$.
Indeed, we have $f'(z) = (1+z^2/2)^{-3/2}>0$ for each $z\in\RR_+$.
This special choice of $f$ plays a role in robust statistics, for more details, see, e.g., Rey \cite[Section 6.4]{Rey} or Example \ref{Ex_Mat}.

In the next remark, restricted to a one-dimensional parameter set, we recall Theorem 3.2 in Clarke \cite{Cla}
 on the local uniqueness of a root of the function \eqref{fgv_elmeleti}.
We will see that the assumptions of Theorem 3.2 in Clarke \cite{Cla} are much more involved and quite different compared to those of our Theorem \ref{Thm_theor_M_est_uniq}.

\begin{Rem}\label{Rem_Clarke}
Let $X:=\RR$, $\Theta$ be a nonempty open interval of $\RR$, and $\psi:\RR\times \Theta\to\RR$ be a measurable function
in its first variable.
Given $t_0\in\Theta$ and a distribution function $G:\RR\to [0,1]$, let
 \[
  I(\psi,G):=\left\{  t\in\Theta : \int_{\RR}\psi(x,t)\,\dd G(x) = 0 \right\},
 \]
 and let $T(\psi,G)\in\Theta$ be a solution of the minimization problem
 \[
   \min_{t\in I(\psi,G)} \vert t-t_0\vert,
 \]
 provided that $I(\psi,G)$ is nonempty.

In statistical estimation theory, a family of distribution functions $\{F_t:\RR\to[0,1]: t\in \Theta\}$ is given, and one chooses $G:=F_{t_0}$ in the minimization problem above.
We also note that the minimization problem above is a special case of a more general one given in (1.3) in Clarke \cite{Cla},
 where a so-called selection function $\varrho: \cE\times \Theta\to\RR$ comes into play,
 where $\cE$ denotes the set of distribution functions on $\RR$.
Namely, the minimization problem in (1.3) in Clarke \cite{Cla} with the selection function $\varrho(G,t):=\vert t-t_0\vert$,
 $G\in\cE$, $t\in\Theta$, gives the minimization problem above.
In Clarke \cite{Cla2} one can find several interesting examples for other selection functionals, for example,
 $\varrho_1(G,t):=\int_{\RR} (G(x)-F_t(x))^2\,\dd K(x)$, $G\in\cE$, $t\in\Theta$, where $K:\RR\to\RR_+$ is a suitable
 weight function, or $\varrho_2(G,t):=\vert \Med(G) - t \vert$, $G\in\cE$, $t\in\Theta$,
 where $\Med(G)$ denotes the median of $G$ (provided that it exists uniquely).
Note that $\varrho_2$ selects the root of $\int_{\RR}\psi(x,t)\,\dd G(x) = 0$, $t\in\Theta$,
 which is the closest to the median of $G$.

Given $t_0\in\Theta$ and a distribution function $G:\RR\to [0,1]$, suppose the following assumptions:
 \begin{enumerate}
    \item[$(A0)$] $t_0\in I(\psi,G)$. In this case, we have $T(\psi,G)=t_0$.
    \item[$(A1)$] $\psi$ has a continuous (partial) derivative with respect to its second variable on $\RR\times D$,
                 where $D\subseteq \Theta$ is a compact interval containing $t_0$ in its interior.
    \item[$(A2)$] there exist a function $g:\RR\to\RR_+$ and $\vare>0$ such that
          \begin{itemize}
             \item $\vert \psi(x,t)\vert \leq g(x)$ for each $x\in \RR$ and $t\in D$,
             \item $\vert \partial_2\psi(x,t)\vert \leq g(x)$ for each $x\in \RR$ and $t\in D$,
             \item $\int_{\RR} g(x)\, \dd F(x) <\infty$ for each $F\in K(G,\vare)$, where
                   $F\in K(G,\vare)$ denotes the open neighbourhood of $G$ with radius $\vare$
                   with respect to the Kolmogorov's distance $d_K$ of distribution functions given by
                   $d_K(F,\widetilde F):=\sup_{x\in\RR} \vert F(x) - \widetilde F(x)\vert$ for
                   distribution functions $F$ and $\widetilde F$.
          \end{itemize}
    \item[$(A3)$] $\int_\RR \partial_2\psi(x,t_0)\,\dd G(x)\ne 0$.
    \item[$(A4)$] for each $\delta>0$ there exists $\vare>0$ such that for each $F\in K(G,\vare)$, we have that
                  \[
                      \sup_{t\in D} \Big\vert  \int_\RR \psi(x,t)\,\dd F(x) - \int_\RR \psi(x,t)\,\dd G(x)  \Big\vert <\delta
                  \]
                  and
                  \[
                      \sup_{t\in D} \Big\vert  \int_\RR \partial_2\psi(x,t)\,\dd F(x) - \int_\RR \partial_2\psi(x,t)\,\dd G(x)  \Big\vert <\delta.
                  \]
                  In particular, Condition (A4) implies that, for each $t\in D$, the functionals $\cE\ni F \mapsto \int_\RR \psi(x,t)\,\dd F(x)$
                   and $\cE\ni F \mapsto \int_\RR \partial_2\psi(x,t)\,\dd F(x)$ are continuous at $G$ with respect to the Kolmogorov's distance $d_K$.
 \end{enumerate}

Given $t_0\in\Theta$ and a distribution function $G:\RR\to [0,1]$, under the assumptions $(A0)-(A4)$,
Clarke \cite[Theorem 3.2]{Cla} proved that for any $\kappa>0$ there exists an $\vare>0$ such that
 $T(\psi,F)$ exists for each $F\in K(G,\vare)$ and $T(\psi,F)\in (t_0-\kappa,t_0+\kappa)$.
Further, for this $\vare$ there exists a $\kappa^*>0$ such that
 \[
   I(\psi,F)\cap (t_0-\kappa^*, t_0+\kappa^*) = T(\psi,F),
 \]
 that is, the equation $\int_{\RR}\psi(x,t)\, F(\dd x)=0$, $t\in\Theta$, has a unique solution in the interval $(t_0-\kappa^*, t_0+\kappa^*)$.
In particular, by choosing $F:=G$, one can see that the equation $\int_{\RR}\psi(x,t)\,\dd G(x)=0$, $t\in \Theta$, has a unique solution {\sl locally}, provided that it has a solution.

We note that, if the Kolmogorov's distance in the assumptions (A2) and (A4) is replaced by the L\'evy's distance for distribution functions, then an analogous statement holds (see Clarke \cite[page 1197]{Cla}).
\proofend
\end{Rem}

\section{Examples from statistical estimation theory}
\label{Section_statistical_examples}

In this section, we present several examples from statistical estimation theory, where our results in Sections \ref{Section_Ex_uniq_empirical} and \ref{Section_Ex_uniq_theoretical} can be well-applied.
For example, we consider the cases of the empirical median, the empirical quantiles, the empirical expectiles, $\psi$-estimators recently used by Mathieu \cite{Mat}, some $\psi$-estimators that are important in robust statistics, and we also study some examples from maximum likelihood theory.
The proof of the results in this section (Propositions \ref{Pro_quantile} and \ref{Pro_Mat})
 can be found in Section \ref{Section_Proof_for_Section4}.

\begin{Ex}[Empirical median]\label{Ex_median}
Let $X:=\RR$, $\Theta:=\RR$ and $\psi:\RR\times\RR\to\RR$, $\psi(x,t):=\sign(x-t)$, $x,t\in\RR$.
For each $x\in \RR$, the function $\RR\ni t\mapsto \psi(x,t)$ is decreasing, but not strictly decreasing.
Then for each $n\in\NN$ and $x_1,\ldots,x_n\in\RR$, the equation \eqref{psi_est_equation} takes the form
 \begin{align}\label{psi_est_equation_median}
   \sum_{i=1}^n \sign(x_i-t) =0,\qquad t\in\RR.
 \end{align}
In this special case, the function $\psi$ is not continuous in its second variable, i.e., $\psi\notin\Psi[Z](\RR,\RR)$, and the corresponding equation \eqref{psi_est_equation_median} has an important role in statistics.
Namely, one can check that \ $\mathrm{Med}_n:\RR^n\to\RR$,
 \begin{align}\label{help_med_2}
     \mathrm{Med}_n(x_1,\ldots,x_n):=
        \frac{1}{2}\left( x^*_{\lceil \frac{n}{2}\rceil} + x^*_{\lfloor \frac{n}{2}+1 \rfloor}\right)
        =\begin{cases}
            x_{k+1}^* & \text{if $n=2k+1$,}\\
            \frac{1}{2}(x_k^* + x_{k+1}^*)  & \text{if $n=2k$,}
        \end{cases}\qquad k\in\ZZ_+,
 \end{align}
is a solution of the equation \eqref{psi_est_equation_median},
where $x_1^*\leq x_2^*\leq\cdots\leq x_n^*$ denotes the ordered sample of $x_1,\ldots,x_n\in\RR$.
Of course, if $n=2k$, then there are other solutions of the equation \eqref{psi_est_equation_median}.
For example, if $x_1<x_2<\cdots < x_{2k}$, then we have
 \begin{align}\label{help_med_1}
   \bigg\{ t\in\RR : \sum_{i=1}^{2k} \sign(x_i-t)   =0 \bigg\}
      = [x_k,x_{k+1}],
 \end{align}
 where $[x_k,x_{k+1}]$ is not a singleton.
Note that $\mathrm{Med}_n(x_1,\ldots,x_n)$ is nothing else but the well-known empirical median of $x_1,\ldots,x_n$.

Further, we have that $\psi\in\Psi[T_n](\RR,\RR)$ for each $n=2k+1$, $k\in\ZZ_+$, with $\vartheta_n(\bx)=x_{k+1}^*$, $\bx=(x_1,\ldots,x_n)\in\RR^n$.
Note also that $\psi\notin\Psi[T_n](\RR,\RR)$ for any $n=2k$, $k\in\NN$.
Indeed, if $x_1<x_2<\cdots < x_{2k}$, then \eqref{help_med_1} implies that
 the function $\psi_{(x_1,\ldots,x_{2k})}$ does not have a point of sign change.
Furthermore, for each $x,y\in\RR$ with $\vartheta_1(x)=x<y=\vartheta_1(y)$, the function \eqref{function_newhanyados} takes the form
 \[
   (x,y)\ni t\mapsto \frac{\psi(x,t)}{\psi(x,t) - \psi(y,t)} = \frac{\sign(x-t)}{\sign(x-t) - \sign(y-t)}
                                                             = \frac{-1}{-1-1} = \frac{1}{2}.
 \]
This function is a rational constant, in particular  strictly $\frac{1}{\ell}$-increasing for each $\ell\in\NN$.
It underlines the fact that the strictly $\frac{1}{n}$-increasing property of the function \eqref{function_newhanyados} in Proposition \ref{Pro_3} is only a necessary, but not a sufficient condition for $\psi\in\Psi[T_n](\RR,\RR)$.
Indeed, in the present example, the function \eqref{function_newhanyados} is strictly $\frac{1}{\ell}$-increasing for each $\ell\in\NN$, but $\psi\notin\Psi[T_n](\RR,\RR)$ for $n=2k$, $k\in\NN$.
Moreover, for each $x,y\in\RR$ with $\vartheta_1(x)=x<y=\vartheta_1(y)$, the function \eqref{function_newhanyados0}  takes the form
 \[
   (x,y)\ni t\mapsto -\frac{\psi(x,t)}{\psi(y,t)} = -\frac{\sign(x-t)}{\sign(y-t)} = 1.
 \]
This function is a rational constant, in particular, increasing, but not strictly increasing.
It is in accordance with part (iv) of Theorem \ref{Thm_M_est_uniq} and part (i) of Proposition \ref{Pro_M_est_uniq} as well,  since $\psi\in\Psi[T_n](\RR,\RR)$ for infinitely many $n\in\NN$, namely, for each $n=2k+1$, $k\in\ZZ_+$;
 and for each $x\in \RR$, the function $\RR\ni t\mapsto \psi(x,t)$ is decreasing.
Note also that it does not hold that for each $m\in\NN$ there exists an $n\in\NN$
 such that $m$ divides $n$ and $\psi\in\Psi[T_n](\RR,\RR)$ (indeed, in case of an even $m\in\NN$ one cannot choose such an $n$).
Moreover, one cannot apply part (vi) of Theorem \ref{Thm_M_est_uniq}.
This underlines that the increasing property of the function \eqref{function_newhanyados0}
 is only a necessary, but not a sufficient condition in order that $\psi$ be a $T_n$-function for each $n\in\NN$.
Finally, we mention that the present example also shows that in Proposition \ref{T_n} the restriction that $m\in\{1,\ldots,n\}$ divides $n$ cannot be removed in general.
\proofend
\end{Ex}

\begin{Ex}[Empirical quantiles]\label{Ex_quantile}
Given $\alpha\in(0,1)$, $n\in\NN$ and $x_1,\ldots,x_n\in\RR$, an empirical $\alpha$-quantile based on $x_1,\ldots,x_n$ is defined as any solution of the minimization problem:
 \begin{align}\label{help_quantile_1}
  \begin{split}
   \min_{t\in\RR} \sum_{i=1}^n \varphi_\alpha(x_i-t)
             &= \min_{t\in\RR} \sum_{i=1}^n \Big(\alpha \bone_{\{x_i\geq t\}} + (\alpha-1) \bone_{\{x_i<t\}} \Big)(x_i-t)\\
             &= \min_{t\in\RR} \sum_{i=1}^n \frac{1}{2}\Big( \vert x_i-t\vert + (2\alpha-1)(x_i-t)\Big),
  \end{split}
 \end{align}
 where \ $\varphi_\alpha:\RR\to\RR$ \ is the so-called $\alpha$-quantile check function given by
 \[
   \varphi_\alpha(x):=\vert \alpha - \bone_{\{x<0\}}\vert \vert x\vert
                     = \big(\alpha \bone_{\{x\geq 0\}} + (\alpha-1)\bone_{\{x<0\}}\big)x
                     = \frac{1}{2}\Big(\vert x\vert + (2\alpha-1)x\Big),\qquad x\in\RR,
 \]
 see, e.g., Koenker and Bassett \cite[Section 3]{KoeBas}.
Some authors call a solution of \eqref{help_quantile_1} as a geometric $\alpha$-quantile, see, e.g.,
 Passeggeri and Reid \cite[Section 2]{PasRei}.

It is known that $q\in\RR$ is an empirical $\alpha$-quantile based on $x_1,\ldots,x_n$ if and only if
 the following two inequalities hold
 \begin{align}\label{help_quantile_2}
   \frac{1}{n}\sum_{i: \, x_i<q} 1 \leq \alpha
      \qquad \text{and}\qquad
    \alpha \leq \frac{1}{n}\sum_{i: \, x_i\leq q } 1,
 \end{align}
 see, e.g., Lange \cite[Problems 12.12/13.]{Lan}.

It is also known that, given $\alpha\in(0,1)$, $n\geq 2$, $n\in\NN$, and $x_1,\ldots,x_n\in\RR$,
an empirical $\alpha$-quantile given as a solution of the minimization problem \eqref{help_quantile_1} is uniquely defined if and only if $\alpha\not\in\{\frac{1}{n}, \frac{2}{n}, \ldots, \frac{n-1}{n}\}$, and in
 case of uniqueness, we have that it is given by
 \begin{align}\label{help_quantile_3}
   q_n^{(\alpha)}(x_1,\ldots,x_n):=\frac{1}{2}(x^*_{\lceil n\alpha \rceil} + x^*_{\lfloor n\alpha+1 \rfloor}),
 \end{align}
where $x_1^*\leq x_2^*\leq\cdots\leq x_n^*$ denotes the ordered sample of $x_1,\ldots,x_n\in\RR$, see Passeggeri and Reid \cite[Lemma 4.1]{PasRei}.
We also mention an interesting result of Passeggeri and Reid \cite[Lemma 4.2]{PasRei}, which states that, given $\alpha\in(0,1)$ and $n\geq 2$, $n\in\NN$, the function $q_n^{(\alpha)}:\RR^n\to\RR$ given by \eqref{help_quantile_3} is Lipschitz continuous:
 \[
   \vert q_n^{(\alpha)}(x_1,\ldots,x_n) - q_n^{(\alpha)}(y_1,\ldots,y_n) \vert \leq \max_{j\in\{1,\ldots,n\}} \vert x_j -y_j\vert
 \]
 for each $x_1,\ldots,x_n,y_1,\ldots,y_n\in\RR$.

Note that if $\alpha=\frac{1}{2}$, then the empirical median $\mathrm{Med}_n(x_1,\ldots,x_n)$ of $x_1,\ldots,x_n$ given in \eqref{help_med_2}
 is a solution of the minimization problem
 \[
   \min_{t\in\RR} \sum_{i=1}^n \vert x_i-t\vert.
 \]
Indeed, if $\alpha=\frac{1}{2}$, then $\varphi_{1/2}(x)=\frac{1}{2}\vert x\vert$, $x\in\RR$, and hence the minimization problem
 \eqref{help_quantile_1} with $\alpha=\frac{1}{2}$ is equivalent to $\min_{t\in\RR} \sum_{i=1}^n \vert x_i-t\vert$, and, as we have recalled, $q\in\RR$ is a solution of this minimization problem if and only if the inequalities \eqref{help_quantile_2} hold for $q$ with $\alpha=\frac{1}{2}$.
Further, one can easily check that the empirical median $\mathrm{Med}_n(x_1,\ldots,x_n)$ based on $x_1,\ldots,x_n$
 satisfies the inequalities \eqref{help_quantile_2} with $\alpha=\frac{1}{2}$.
Moreover, for $\alpha=\frac{1}{2}$ and $n\geq 2$, $n\in\NN$, we have that
 $\alpha\not\in\{\frac{1}{n}, \frac{2}{n}, \ldots, \frac{n-1}{n}\}$ holds if and only if
 $n=2k+1$ with some $k\in\NN$, and in this case
 \[
  q_n^{(\frac{1}{2})}(x_1,\ldots,x_n) = \frac{1}{2}(x^*_{\lceil \frac{2k+1}{2} \rceil} + x^*_{\lfloor \frac{2k+1}{2} +1 \rfloor})
        = \frac{1}{2}(x^*_{k+1} + x^*_{k+1})
        = x^*_{k+1}
        = \mathrm{Med}_n(x_1,\ldots,x_n)
 \]
 for all $x_1,\ldots,x_n\in\RR$, where $\mathrm{Med}_n(x_1,\ldots,x_n)$ is defined in \eqref{help_med_2}.
\proofend
\end{Ex}

Motivated by the minimization problem \eqref{help_quantile_1}, we investigate the function $\psi:\RR\times\RR\to\RR$,
 \begin{align}\label{help17}
   \psi(x,t):=\begin{cases}
                \alpha & \text{if $x>t$,}\\
                0 & \text{if $x=t$,}\\
                \alpha -1 & \text{if $x<t$.}
               \end{cases}
 \end{align}
For each $x\in \RR$, the function $\RR\ni t\mapsto \psi(x,t)$ is decreasing, but not strictly decreasing, and not continuous.
By choosing $X:=\Theta:=\RR$, we have that $\psi$ is a $T_1$-function with $\vartheta_1(x):=x$, $x\in\RR$.
Analogously to the result of Passeggeri and Reid \cite[Lemma 4.1]{PasRei}, we are going to show the following result.

\begin{Pro}\label{Pro_quantile}
Given $\alpha\in(0,1)$, for each $n\geq 2$, the function $\psi$ defined by \eqref{help17} has the property $[T_n]$ if and only if $\alpha\not\in\big\{\frac{1}{n},\dots,\frac{n-1}{n}\big\}$.
Further, in this case, for all $x_1,\dots,x_n\in\RR$, we have that
\[
 \vartheta_n(x_1,\ldots,x_n)
 =x^*_{\lceil n\alpha \rceil}=x^*_{\lfloor n\alpha+1 \rfloor},
\]
and hence \eqref{help_quantile_3} is also valid,
where $x_1^*\leq x_2^*\leq\cdots\leq x_n^*$ denotes the ordered sample of $x_1,\ldots,x_n$.
\end{Pro}

\begin{Ex}[Expectiles]\label{Ex_expectile}
Let $\alpha\in(0,1)$, $n\in\NN$ and $x_1,\ldots,x_n\in\RR$.
The empirical $\alpha$-expectile based on $x_1,\ldots,x_n$ is defined as any solution of the minimization problem:
 \[
   \min_{t\in\RR} \sum_{i=1}^n \widetilde\varphi_\alpha(x_i-t)
             = \min_{t\in\RR} \sum_{i=1}^n \Big(\alpha \bone_{\{x_i\geq t\}} + (1-\alpha) \bone_{\{x_i<t\}} \Big)(x_i-t)^2,
 \]
 where \ $\widetilde\varphi_\alpha:\RR\to\RR$ \ is given by
 \[
   \widetilde \varphi_\alpha(x):=\vert \alpha - \bone_{\{x<0\}}\vert x^2
                     = \big(\alpha \bone_{\{x\geq 0\}} + (1-\alpha)\bone_{\{x<0\}}\big)x^2,\qquad x\in\RR,
 \]
 see, e.g., Newey and Powell \cite{NewPow}.
Expectiles are also called smoothed versions of quantiles or least asymptotically weighted squares estimators.

Motivated by this minimization problem, we may investigate the applicability of Theorem \ref{Thm_M_est_uniq} for the function $\psi:\RR\times\RR\to\RR$,
 \begin{equation}\label{Exp-theor}
   \psi(x,t):=\alpha (x-t)^+-(1-\alpha)(x-t)^-
   =\begin{cases}
                \alpha(x-t) & \text{if $x>t$,}\\
                0 & \text{if $x=t$,}\\
                (1-\alpha)(x-t) & \text{if $x<t$.}
               \end{cases}
 \end{equation}
For each $x\in \RR$, the function $\RR\ni t\mapsto \psi(x,t)$ is strictly decreasing.
By choosing $X:=\RR$ and $\Theta:=\RR$, we have $\psi$ is a $T_1$-function with $\vartheta_1(x):=x$, $x\in\RR$,
 and for each $x,y\in\RR$ with $\vartheta_1(x)=x<y=\vartheta_1(y)$, the function \eqref{function_newhanyados0}  takes the form
 \[
   (x,y)\ni t\mapsto -\frac{\psi(x,t)}{\psi(y,t)} = -\frac{(1-\alpha)(x-t)}{\alpha(y-t)}
                                                 = - \frac{1-\alpha}{\alpha}\left(1-\frac{y-x}{y-t}\right),
 \]
which is a strictly increasing function.
It is in accordance with part (i) of Proposition \ref{Pro_M_est_uniq}.
Since $\psi(x,\vartheta_1(x))=0$, $x\in\RR$, by part (vi) of Theorem \ref{Thm_M_est_uniq}, we have $\psi\in\Psi[T_n^\blambda](X,\Theta)$ for each $n\in\NN$ and $\blambda\in\Lambda_n$.
In particular, using that $\psi$ is continuous in its second variable, we also have, for each $n\in\NN$ and $x_1,\ldots, x_n\in\RR$, that the equation $\sum_{i=1}^n \psi(x_i,t) = 0$, $t\in\RR$,
has a unique solution, see part (ii) of Remark \ref{Rem_1}.
\proofend
\end{Ex}

Now, using Proposition \ref{Pro_theor_M_est_uniq2}, we show that the $\psi$-expectation function with the given $\psi$ and a random variable
 $\xi$ such that $\EE(\vert \xi\vert)<\infty$ has a unique zero, which is known to be the $\alpha$-expectile of $\xi$,
 see, e.g., Bellini et al.\ \cite[Section 2]{BelKlaMul}.

\begin{Pro}\label{Pro_expectile}
Let $\psi\in\Psi(\RR,\RR)$ be defined by \eqref{Exp-theor}, let $\alpha\in(0,1)$ and $\xi$ be a random variable on a probability space $(\Omega,\cA,\PP)$ such that $\EE(\vert \xi\vert)<\infty$.
Then the equation $\EE(\psi(\xi,t))=0$, $t\in\RR$, has a unique solution, which is known to be the $\alpha$-expectile of $\xi$.
\end{Pro}

We note that Proposition \ref{Pro_expectile} also follows from Lemma A.1 in Kr\"atschmer and Z\"ahle \cite{KraZah},
 where it was shown that the mapping $\RR\ni t\mapsto \EE(\psi(\xi,t))$ is real-valued, continuous, strictly decreasing and it satisfies that
 $\lim_{t\to\pm\infty} \EE(\psi(\xi,t)) = \mp\infty$.
Nonetheless, we give a proof of Proposition \ref{Pro_expectile} using Proposition \ref{Pro_theor_M_est_uniq2}, since we would like to
 demonstrate the applicability of our result.
Note also that if $n\in\NN$ and $x_1,\ldots,x_n\in\RR$, then, by choosing $\xi$ as a random variable such that
 $\PP(\xi=x_i)=\frac{1}{n}$, $i\in\{1,\ldots,n\}$, then $\EE(\psi(\xi,t))=\frac{1}{n}\sum_{i=1}^n \psi(x_i,t)$, $t\in\RR$,
 so Proposition \ref{Pro_expectile} yields that the equation $\sum_{i=1}^n \psi(x_i,t)=0$, $t\in\RR$, has a unique solution.

Next, we recall a function $\psi$ that has been recently used
for constructing M-estimators by Mathieu \cite{Mat} (see also Example \ref{Ex_Mat_elmeleti}), and, in case of the Huber function, we take advantage of Proposition \ref{T_n}.

\begin{Ex}\label{Ex_Mat}
Let $X:=\RR$, $\Theta:=\RR$ and $\psi:\RR\times\RR\to\RR$,
 \begin{align}\label{help_psi_mathieu}
   \psi(x,t):=\sign(x-t)f(\vert x-t\vert),\qquad x,t\in\RR,
 \end{align}
where $f:\RR_+\to\RR$.
Note that the function $\psi$ given in \eqref{help_psi_mathieu} has the same form as the function given in \eqref{help_psi_mathieu_new},
 the only difference is that the function $f$ can take negative values in case of \eqref{help_psi_mathieu},
 but not in case of  \eqref{help_psi_mathieu_new}.
Then $\psi(x,t) = f((x-t)^+) - f((t-x)^+)$, $x,t\in\RR$.
Note that if $f$ is continuous and $f(0)=0$, then $\psi$ is continuous in its second variable, i.e., $\psi\in\Psi[C](\RR,\RR)$.
We now recall some known special choices for the function $f$ appearing in \eqref{help_psi_mathieu} such that the corresponding $\psi$-estimator has an important role in (robust) statistics:
 \begin{itemize}
  \item[(i)] the Huber function $f_H:\RR_+\to\RR$,
             \begin{align}\label{func_Huber}
              f_H(z):= z\bone_{\{z\leq \beta\}} + \beta \bone_{\{z>\beta\}},\qquad z\in\RR_+,
             \end{align}
             where $\beta>0$ (see Huber \cite{Hub64}), which is a continuous and increasing (but not strictly increasing) function starting from $0$.
             Then the function $\psi:\RR\times\RR\to\RR$ given in \eqref{help_psi_mathieu} takes the form
             \[
              \psi(x,t)= \begin{cases}
                          x-t & \text{if $\vert x-t\vert\leq \beta$,}\\
                          \beta\sign(x-t) & \text{if $\vert x-t\vert>\beta$,}
                        \end{cases} \qquad x,t\in\RR.
             \]
              In general, $\psi$ is not a $T_2$-function.
              Indeed, for example, if $\beta:=1$, $x_1:=0$ and $x_2:=3$, then
             \[
              \psi_{(x_1,x_2)}(t)=
                \sum_{i=1}^2 \psi(x_i,t)
               = \begin{cases}
                  1+1=2>0 & \text{if $t\leq -1$,}\\
                  -t+1>0 & \text{if $-1<t\leq 1$,}\\
                  -1+1=0 & \text{if $1<t\leq 2$,}\\
                  -1+3-t=2-t<0 & \text{if $2<t<4$,}\\
                 -1-1=-2<0 & \text{if $t\geq 4$.}
                \end{cases}
 \]
Consequently, the function $\psi_{(x_1,x_2)}$ does not have a point of sign change, and, thus $\psi$ is not a $T_2$-function,
 and, by Proposition \ref{T_n}, it yields that $\psi$ is not a $T_{2k}$-function for any $k\in\NN$.
On the other hand, for each $k\in\NN$, $\psi$ is a $T_{2k+1}$-function.
Indeed, let $k\in\NN$ and $x_1,\dots,x_{2k+1}\in\RR$ be arbitrary.
Then the function $\RR\ni t\mapsto \sum_{i=1}^{2k+1} \psi(x_i,t)$ is decreasing (because each of its terms is decreasing), and, on the contrary, let us assume that it does not have a point of sign change.
Then it is constant on a proper subinterval $I$ of $\RR$.
In this case, each term in question must be also constant on this subinterval, and, by taking into account the form of $\psi$,
it must be equal to $\beta$ or to $-\beta$ on $I$.
However, a $(2k+1)$-term sum whose terms are equal to $\beta$ or to $-\beta$ cannot be equal to zero, which leads us to a contradiction.
All in all, $\psi$ is a $T_n$-function if $n\in\NN$ is odd, and $\psi$ is not a $T_n$-function if $n\in\NN$ is even.

 \item[(ii)] the Catoni function $f_C:\RR_+\to\RR$,
             \[
              f_C(z):=\ln\left(1+\frac{z}{b}+ \frac{1}{2} \left(\frac{z}{b}\right)^2 \right),\qquad z\in\RR_+,
             \]
             where $b>0$ (see Catoni \cite{Cat12}), which is continuous and strictly increasing starting from $0$.
 \item[(iii)] a polynomial function $f_P:\RR_+\to\RR$,
             \[
               f_P(z):=\frac{z}{1+\big(\frac{z}{\beta}\big)^{1-\frac{1}{p}}},\qquad z\in\RR_+,
             \]
             where $p\in\NN$ and $\beta>0$, which is continuous and strictly increasing starting from $0$.
 \item[(iv)] another Catoni-type function $\widetilde f_C:\RR_+\to\RR$,
             \[
               \widetilde f_C(z):= \ln\left(1+z+\frac{z^\alpha}{\alpha} \right),\qquad z\in\RR_+,
             \]
             where $\alpha\in(1,2)$ (see Chen et al.\ \cite{CheJinLiXu}), which is continuous and strictly increasing starting from $0$.
  \item[(v)] $f:\RR_+\to\RR$, $f(z):=z/\sqrt{1+z^2/2}$, $z\in\RR_+$, which is a continuous and strictly increasing function starting from $0$.
             Indeed, we have $f'(z) = (1+z^2/2)^{-3/2}>0$ for each $z\in\RR_+$.
             Then the function $\psi:\RR\times\RR\to\RR$ given in \eqref{help_psi_mathieu}
             takes the form
             \[
               \psi(x,t)= \frac{x-t}{\sqrt{1+\frac{(x-t)^2}{2}}}, \qquad x,t\in\RR.
             \]
             In robust statistics, one calls $\psi$ as the $L_1$-$L_2$ function, see, e.g., Rey \cite[Section 6.4]{Rey}.
  \item[(vi)] $f:\RR_+\to\RR$, $f(z):=z/(1+z)$, $z\in\RR_+$, which is a continuous and strictly increasing function starting from $0$.
              Then the function $\psi:\RR\times\RR\to\RR$ given in \eqref{help_psi_mathieu} takes the form
              \[
               \psi(x,t)= \frac{x-t}{1+\vert x-t\vert}, \qquad x,t\in\RR,
               \]
              which is called a ''fair''-type function in robust statistics, see, e.g., Rey \cite[Section 6.4]{Rey}.
 \end{itemize}
\proofend
\end{Ex}

Next, we discuss the applicability of Theorem \ref{Thm_M_est_uniq} for the function $\psi$ given in \eqref{help_psi_mathieu}.
Namely, we prove the following statement.

\begin{Pro}\label{Pro_Mat}
Let $f:\RR_+\to\RR$ be a function with $f(0)=0$ and let $\psi$ be given by \eqref{help_psi_mathieu}.
Then we have that
\vspace*{-3mm}
\begin{itemize}
 \item[(a)] $\psi\in\Psi[T_1](\RR,\RR)$ if and only if $f(z)>0$ for all $z>0$, and, in this case, $\vartheta_1(x)=x$ and $\psi(x,\vartheta_1(x))=0$ hold for all $x\in\RR$.
 \item[(b)] if $\psi\in\Psi[T_n](\RR,\RR)$ for infinitely many $n\in\NN$, then $f$ is increasing.
 \item[(c)] $\psi\in\Psi[T_2^{\blambda}](\RR,\RR)$ for each $\blambda\in\Lambda_2$ if and only if $f$ is strictly increasing.
 \item[(d)] $\psi\in\Psi[T_n^{\blambda}](\RR,\RR)$ for each $n\in\NN$ and $\blambda\in\Lambda_n$ if and only if $f$ is strictly increasing.
\end{itemize}
\end{Pro}

As a consequence of Proposition \ref{Pro_Mat}, if $f:\RR_+\to\RR$ is a strictly increasing and continuous function such that $f(0)=0$, then $\psi$ given in \eqref{help_psi_mathieu} is continuous in its second variable, and hence for each $n\in\NN$ and $x_1,\ldots,x_n\in \RR$, the equation $\sum_{i=1}^n \psi(x_i,t) = 0$, $t\in \RR$ (i.e., the equation \eqref{psi_est_equation})
 has a unique solution.

In the special cases (ii)--(vi) of Example \ref{Ex_Mat}, by part (d) of Proposition \ref{Pro_Mat},
 the corresponding function $\psi$ given in \eqref{help_psi_mathieu} is a $T_n^\blambda$-function for each $n\in\NN$ and $\blambda\in\Lambda_n$.
In particular, since in the cases (ii)-(vi), $\psi$ is continuous in its second variable,
 for each $n\in\NN$ and $x_1,\ldots,x_n\in\RR$, the equation \eqref{psi_est_equation} with the given function $\psi$ has a unique solution, see part (ii) of Remark \ref{Rem_1}.

In the remaining part of this section, we investigate the applicability of Theorem \ref{Thm_M_est_uniq} for finding solutions of likelihood equations in the theory of MLEs.
Let $\Theta$ be a nonempty open interval of $\RR$, and $f:\RR\times \Theta\to\RR$ be a function such that for each $t\in\Theta$, \ the function $\RR\ni x \mapsto f(x,t)$ is a density function.
Let us introduce the set
 \[
 \cX_f:=\{x\in\RR:f(x,t)>0,\ \forall\;t\in\Theta\},
 \]
 and suppose that $\cX_f$ is nonempty.
Note that, in general, it can happen that $\cX_f=\emptyset$.
For example, if $\Theta=(0,\infty)$ and $f:\RR\times (0,\infty)\to\RR$,
 \[
   f(x,t):=\begin{cases}
              \frac{1}{t} & \text{if $x\in(0,t)$},\\
              0           &  \text{if $x\not\in(0,t)$},
           \end{cases} \qquad t\in(0,\infty),
 \]
i.e, for each $t\in(0,\infty)$, the function $\RR\ni x\mapsto f(x,t)$ is the density function of a uniformly distributed random variable on $(0,t)$, then $\cX_f=\emptyset$.
Turning back to the case when $\cX_f\ne\emptyset$,
and supposing that the (partial) derivative $\partial_2 f$ of $f$ with respect to its second variable exists,
the equation \eqref{psi_est_equation} with $X:=\cX_f$ and the function $\psi:\cX_f\times\Theta\to\RR$ defined by
 \begin{align}\label{psi_MLE}
   \psi(x,t):=\partial_2(\ln(f(x,t)))=\frac{\partial_2 f(x,t)}{f(x,t)},   \qquad (x,t)\in \cX_f\times\Theta,
 \end{align}
is nothing else but the likelihood equation based on the observations $x_1,\ldots,x_n\in \cX_f$ in the theory of MLEs.
In some cases, we need to consider an appropriate Borel subset $\widetilde\cX_f$ of $\cX_f$ such that
 $\PP(\xi_t \in \widetilde\cX_f)=1$ for all $t\in\Theta$, where $\xi_t$ is a random variable having a density function $\RR\ni x\mapsto f(x,t)$.
For such a case, see the second part of Example \ref{Ex_1_norm}.

In the next examples, we demonstrate the applicability of Theorem \ref{Thm_M_est_uniq} together
 with Proposition \ref{Pro_M_est_uniq} for proving existence and uniqueness of a solution of the likelihood equation \eqref{psi_est_equation} corresponding to the function $\psi$ given in \eqref{psi_MLE}.

\begin{Ex}\label{Ex_1_norm}
Let $\xi$ be a normally distributed random variable with mean $m\in\RR$ and with variance $\sigma^2$, where $\sigma>0$.
Let $n\in\NN$ and $x_1,\ldots,x_n\in\RR$ be a realization of a sample of size $n$ for $\xi$.
Here by a sample of size $n$, we mean independent and identically distributed random variables \ $\xi_1,\ldots,\xi_n$ \
 with common distribution as that of $\xi$.
It is known that, supposing that $\sigma$ is known, there exists a unique MLE of $m$ based on $x_1,\ldots,x_n\in\RR$,
 and it takes the form \ $\widehat m_n:=\frac{x_1+\cdots+x_n}{n}$.
We will establish the existence and uniqueness of a solution of the corresponding likelihood equation using Theorem \ref{Thm_M_est_uniq} together
 with Proposition \ref{Pro_M_est_uniq}.
In this case, we have $\Theta=\RR$ and $f:\RR\times\RR\to\RR$,
 \[
  f(x,m):=\frac{1}{\sqrt{2\pi}\sigma} \ee^{-\frac{(x-m)^2}{2\sigma^2}}, \qquad x, m\in\RR,
 \]
and consequently $\cX_f=\RR$.
Then $\psi:\RR\times\RR\to\RR$,
 \[
    \psi(x,m)=\frac{\frac{1}{\sigma^2}(x-m)f(x,m)}{f(x,m)} = \frac{1}{\sigma^2}(x-m),\qquad x,m\in\RR.
 \]
Hence $\psi\in\Psi[C,Z_1](\RR,\RR)$ with $\vartheta_1(x):=x$, $x\in\RR$, and $\psi$ is strictly decreasing in its second variable.
Further, using Proposition \ref{Pro_M_est_uniq} and Theorem \ref{Thm_M_est_uniq} (with $X:=\cX_f=\RR$), we can conclude that
 for each $n\in\NN$ and $x_1,\ldots,x_n\in \RR$, the (likelihood) equation \eqref{psi_est_equation}
 with the given function $\psi$ has a unique solution, which is equal to
 $\vartheta_n(x_1,\dots,x_n)=\frac{x_1+\cdots+x_n}{n}=\widehat m_n$, as desired.

It is also known that, supposing that $m$ is known, there exists a unique MLE of $\sigma^2$ based on $x_1,\ldots,x_n\in \RR\setminus\{m\}$,
 and it takes the form \ $\widehat{\sigma^2_n}:=\frac{1}{n}\sum_{i=1}^n (x_i-m)^2$.
We will establish the existence and uniqueness of a solution of the corresponding likelihood equation using Theorem \ref{Thm_M_est_uniq}.
In this case, we have $\Theta=(0,\infty)$, and $f:\RR\times(0,\infty)\to\RR$,
 \[
  f(x,\sigma^2):=\frac{1}{\sqrt{2\pi\sigma^2}} \ee^{-\frac{(x-m)^2}{2\sigma^2}}, \qquad x\in\RR,\; \sigma^2>0,
 \]
and consequently $\cX_f=\RR$.
Instead of $\cX_f$, let us consider its subset $\widetilde\cX_f:=\cX_f\setminus\{m\} = \RR\setminus\{m\}$.
Then $\PP(\xi_{\sigma^2}\in \widetilde\cX_f)=1$ for all $ \sigma^2>0$, where $\xi_{\sigma^2}$ is a normally distributed random variable with mean $m$ and variance $\sigma^2$.
Then $\psi:\widetilde\cX_f\times(0,\infty)\to\RR$,
 \begin{align*}
    \psi(x,\sigma^2)
     & = \frac{1}{f(x,\sigma^2)}\left( \frac{1}{\sqrt{2\pi}}\left(-\frac{1}{2}\right)(\sigma^2)^{-\frac{3}{2}} \ee^{-\frac{(x-m)^2}{2\sigma^2}}
                + \frac{1}{\sqrt{2\pi\sigma^2}} \ee^{-\frac{(x-m)^2}{2\sigma^2}} \frac{(x-m)^2}{2} (\sigma^2)^{-2}
                               \right) \\
     & = \frac{1}{2(\sigma^2)^2}\left( (x-m)^2 - \sigma^2\right) ,
         \qquad x\in\widetilde\cX_f,\; \sigma^2>0.
 \end{align*}
Hence $\psi\in\Psi[C,Z_1](\widetilde\cX_f,(0,\infty))$ with $\vartheta_1(x):=(x-m)^2$, $x\in\widetilde\cX_f$.
Note that $\vartheta_1(x)\in \Theta=(0,\infty)$ for all $x\in\widetilde\cX_f$. (This explains the restriction of $\cX_f=\RR$ to $\widetilde\cX_f=\RR\setminus\{m\}$.)
Further, for each $x,y\in\widetilde\cX_f$ with $\vartheta_1(x)< \vartheta_1(y)$, i.e., $(x-m)^2 < (y-m)^2$, we get that the function
 \begin{align*}
  \big( (x-m)^2, (y-m)^2\big)\ni\sigma^2
     \mapsto -\frac{\psi(x,\sigma^2)}{\psi(y,\sigma^2)}
              = - \frac{(x-m)^2-\sigma^2}{(y-m)^2-\sigma^2}
              = -1 + \frac{(y-m)^2 - (x-m)^2 }{(y-m)^2 - \sigma^2}
 \end{align*}
 is strictly increasing.
Consequently, by part (vi) of Theorem \ref{Thm_M_est_uniq} (with $X:=\widetilde\cX_f$), we conclude that for each $n\in\NN$ and $x_1,\ldots,x_n\in \widetilde\cX_f$, the (likelihood) equation \eqref{psi_est_equation} with the given $\psi$ has a unique solution,
 which is equal to $\vartheta_n(x_1,\dots,x_n)=\frac{1}{n}\sum_{i=1}^n (x_i-m)^2=\widehat{\sigma^2_n}$, as desired.
Finally, we present an alternative argument.
Note that the equation \eqref{psi_est_equation} with the given function $\psi$ has a solution if and only if the equation \eqref{psi_est_equation} with the function
 $\widetilde\psi: \widetilde\cX_f\times(0,\infty)\to\RR$,
 \[
  \widetilde\psi(x,\sigma^2):= 2(\sigma^2)^2 \psi(x,\sigma^2)
                             =(x-m)^2 - \sigma^2,
                             \qquad x\in\widetilde\cX_f, \; \sigma^2>0,
 \]
has a solution, and the two sets of solutions coincide.
Further,  $\widetilde\psi$ is a $T_1$-function, and, for each $x\in\widetilde\cX_f$, the unique point of sign change of the function $(0,\infty)\ni\sigma^2 \mapsto \widetilde\psi(x,\sigma^2)$ is equal to $(x-m)^2$.
The function $\widetilde\psi$ is strictly decreasing in its second variable, and hence part (i) of Proposition \ref{Pro_M_est_uniq} can be applied to $\widetilde\psi$.
This, together with part (vi) of Theorem \ref{Thm_M_est_uniq}, yield that $\widetilde\psi\in\Psi[T_n](\widetilde\cX_f,(0,\infty))$
 for each $n\in\NN$, and hence, trivially, $\psi\in\Psi[T_n](\widetilde\cX_f,(0,\infty))$ for each $n\in\NN$ as well, as expected.
\proofend
\end{Ex}

\begin{Ex}\label{Ex_1_ism}
Let $\alpha>0$ and let $\xi$ be an absolutely continuous random variable with a density function
 \[
   f_\xi(x):=\begin{cases}
               2\alpha x (1-x^2)^{\alpha-1}  & \text{if $x\in(0,1)$,}\\
               0 & \text{otherwise.}
              \end{cases}
 \]
Then one can check that given $n\in\NN$ and a realization $x_1,\ldots,x_n\in(0,1)$
 of a sample of size $n$ for $\xi$, there exists a unique MLE of $\alpha$ and it takes the form
 \[
   \widehat\alpha_n:=-\frac{n}{\sum_{i=1}^n \ln(1-x_i^2)}.
 \]
We will establish the existence and uniqueness of a solution of the corresponding likelihood equation using Theorem \ref{Thm_M_est_uniq} together with Proposition \ref{Pro_M_est_uniq}.
In this case, we have $\Theta=(0,\infty)$  and $f:\RR\times(0,\infty)\to\RR$,
 \[
  f(x,\alpha):=\begin{cases}
               2\alpha x (1-x^2)^{\alpha-1}  & \text{if $x\in(0,1)$, $\alpha>0$,}\\
               0 & \text{otherwise,}
              \end{cases}
 \]
and consequently $\cX_f=(0,1)$.
Then $\psi:(0,1)\times(0,\infty)\to\RR$,
 \[
    \psi(x,\alpha)=\frac{2x\Big((1-x^2)^{\alpha-1} + \alpha (1-x^2)^{\alpha-1}\ln(1-x^2) \Big)}
                    {2\alpha x (1-x^2)^{\alpha-1}}
             = \frac{1}{\alpha} + \ln(1-x^2), \quad x\in(0,1),\; \alpha>0.
 \]
Hence $\psi\in\Psi[C,Z_1]((0,1),(0,\infty))$ with $\vartheta_1(x):=-\frac{1}{\ln(1-x^2)}$, $x\in(0,1)$, and $\psi$ is strictly decreasing in its second variable.
Further, using Theorem \ref{Thm_M_est_uniq} and Proposition \ref{Pro_M_est_uniq} (with $X:=\cX_f=(0,1)$), we can conclude that
 for each $n\in\NN$ and $x_1,\ldots,x_n\in (0,1)$, the (likelihood) equation \eqref{psi_est_equation}
 with the given \ $\psi$ \ has a unique solution,
 which is equal to $\vartheta_n(x_1,\dots,x_n)=-\frac{n}{\sum_{i=1}^n \ln(1-x_i^2)}=\widehat\alpha_n$, as desired.
\proofend
\end{Ex}

\begin{Ex}\label{Ex_1_norm_kevert}
Let $\xi$ be an absolutely continuous random variable with a density function
 \[
   f_\xi(x):=
          \frac{1}{2}\cdot \frac{1}{\sqrt{2\pi}} \ee^{-\frac{x^2}{2}}
                     + \frac{1}{2}\cdot \frac{1}{\sqrt{2\pi}\sigma} \ee^{-\frac{(x-m)^2}{2\sigma^2}} ,
                     \qquad x\in\RR,
 \]
 where $m\in\RR$ and $\sigma>0$.
Note that $f_\xi$ is a mixture density function of the standard normal density function and
 the density function of a normally distributed random variable with mean $m$ and variance $\sigma^2$ with equal $\frac{1}{2}$ weights.
Let $n\in\NN$ and $x_1,\ldots,x_n\in\RR$ \ be a realization of a sample of size $n$ for $\xi$.
In what follows, we assume that $\sigma$ is known, and, using Theorem \ref{Thm_M_est_uniq}, we show that, in general, the corresponding likelihood equation for $m$ may have more solutions.
In this case, we have $\Theta=\RR$ and $f:\RR\times\RR\to\RR$,
 \[
  f(x,m):=\frac{1}{2}\cdot \frac{1}{\sqrt{2\pi}} \ee^{-\frac{x^2}{2}}
                     + \frac{1}{2}\cdot \frac{1}{\sqrt{2\pi}\sigma} \ee^{-\frac{(x-m)^2}{2\sigma^2}}, \qquad x, m\in\RR,
 \]
and consequently $\cX_f=\RR$.
Then $\psi:\RR\times\RR\to\RR$,
 \[
    \psi(x,m)= \frac{\frac{1}{2\sigma^3\sqrt{2\pi}}(x-m)\ee^{-\frac{(x-m)^2}{2\sigma^2}}}{f(x,m)},
                     \qquad x,m\in\RR.
 \]
Hence $\psi\in\Psi[C,T_1](\RR,\RR)$ with $\vartheta_1(x):=x$, $x\in\RR$.
Consequently, for each $n\in\NN$ and $\bx=(x_1,\ldots,x_n)\in X^n$, the likelihood equation \eqref{psi_est_equation} has at least one solution in $\RR$.
Indeed, if $m<\min(x_1,\ldots,x_n)$, then $\sum_{i=1}^n \psi(x_i,m)>0$, and if $m>\max(x_1,\ldots,x_n)$, then $\sum_{i=1}^n \psi(x_i,m)<0$, and hence the continuity of $\psi$ in its second variable together with the Bolzano theorem imply the existence of a solution of \eqref{psi_est_equation}, as desired.
Further, for each $x\in\RR$, we have $\lim_{m\to\pm\infty} \psi(x,m)=0$.
In what follows, we check that it is not true that for each $x,y\in\RR$ with $\vartheta_1(x)<\vartheta_1(y)$, i.e., $x<y$, the function \eqref{function_newhanyados0}
 is increasing.
For each $x,y\in\RR$ with $x<y$, the function \eqref{function_newhanyados0} takes the form
 \begin{align*}
 (x,y)\ni m\mapsto - \frac{\psi(x,m)}{\psi(y,m)}
     = \frac{m-x}{y-m}\cdot \frac{f(y,m) }{f(x,m)} \cdot \frac{ \ee^{\frac{(y-m)^2}{2\sigma^2}}  }{\ee^{\frac{(x-m)^2}{2\sigma^2}} }
     = \frac{(m-x)\Big(\sigma \ee^{-\frac{y^2}{2} + \frac{(y-m)^2}{2\sigma^2} } +1\Big) }{(y-m) \Big(\sigma \ee^{-\frac{x^2}{2} + \frac{(x-m)^2}{2\sigma^2} } +1\Big) },
 \end{align*}
 and, for each $m\in(x,y)$, one can check that
 \begin{align}\label{help15}
  \begin{split}
  \frac{\dd}{\dd m}\left( - \frac{\psi(x,m)}{\psi(y,m)} \right)
     &= \frac{1}{(y-m)^2 \Big(\sigma \ee^{-\frac{x^2}{2} + \frac{(x-m)^2}{2\sigma^2} } +1\Big)^2 }\\
     &\phantom{=\;} \times \Bigg[  (y-x) \Big(\sigma \ee^{-\frac{x^2}{2} + \frac{(x-m)^2}{2\sigma^2} } +1\Big)\Big(\sigma \ee^{-\frac{y^2}{2} + \frac{(y-m)^2}{2\sigma^2} } +1\Big) \\
     &\phantom{=\;\times \Bigg[} -\frac{1}{\sigma} (m-x)(y-m)^2 \ee^{-\frac{y^2}{2} + \frac{(y-m)^2}{2\sigma^2} } \Big(\sigma \ee^{-\frac{x^2}{2} + \frac{(x-m)^2}{2\sigma^2} } +1\Big) \\
      &\phantom{=\;\times \Bigg[} -\frac{1}{\sigma} (m-x)^2(y-m) \ee^{-\frac{x^2}{2} + \frac{(x-m)^2}{2\sigma^2} } \Big(\sigma \ee^{-\frac{y^2}{2} + \frac{(y-m)^2}{2\sigma^2} } +1\Big)
                          \Bigg].
   \end{split}
 \end{align}
 If for each $x,y\in\RR$ with $x<y$, the function \eqref{function_newhanyados0} were increasing, then we would have that
 \[
   \frac{\dd}{\dd m}\left( - \frac{\psi(x,m)}{\psi(y,m)} \right) \geq 0, \qquad m\in(x,y),\quad x<y,\; x,y\in\RR,
 \]
 which, by \eqref{help15}, is equivalent to
 \[
   \sigma(y-x) \geq
   \frac{(m-x)(y-m)^2}{\sigma +\ee^{\frac{y^2}{2} -\frac{(y-m)^2}{2\sigma^2} }} + \frac{(m-x)^2(y-m)}{\sigma +\ee^{\frac{x^2}{2} -\frac{(x-m)^2}{2\sigma^2} }},
    \qquad m\in(x,y), \quad x<y,\; x,y\in\RR.
 \]
However, this inequality does not hold in general, since its left hand side tends to $0$ as $\sigma\downarrow 0$,
 but its right hand side tends to $\infty$ as $\sigma\downarrow 0$.
To give an example, for example, on Figure \ref{Fig_mixed_norm}, we plotted the function
 $(x,y)\ni m\mapsto - \frac{\psi(x,m)}{\psi(y,m)}$ with $x=1$, $y=5$ and $\sigma=1$, which is not increasing.
\begin{figure}[ht!]
 \centering
 \includegraphics[width=10cm,height=8cm]{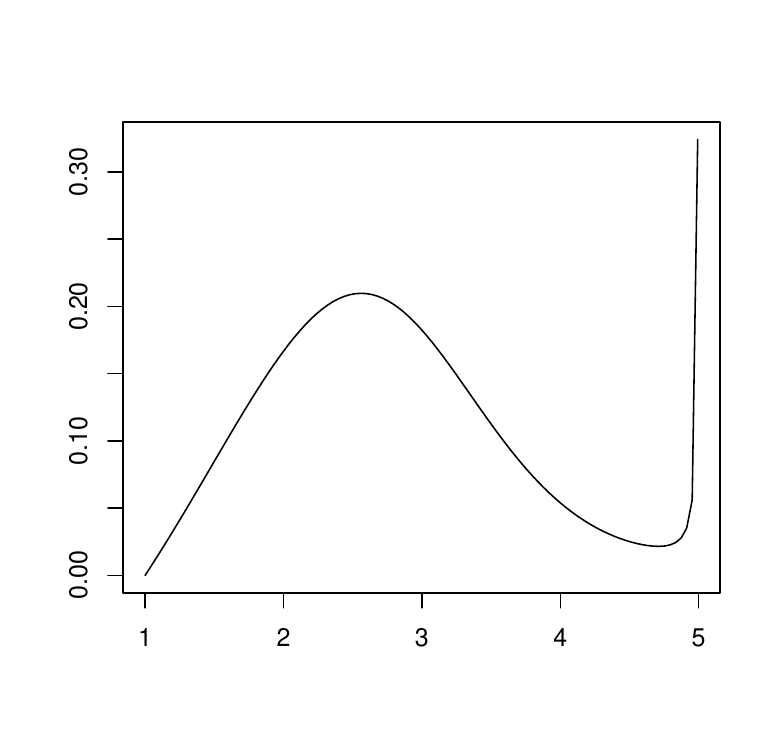}
  \caption{The function $(1,5)\ni m\mapsto - \frac{\psi(x,m)}{\psi(y,m)}$ with $\sigma=1$.}
 \label{Fig_mixed_norm}
\end{figure}
In general, it is not true that for each $x,y\in\RR$ with $x<y$, the function \eqref{function_newhanyados0} is increasing.
If the function \eqref{function_newhanyados0} is not increasing for some $x,y\in\RR$ with $x<y$, then, by part (iv) of Theorem \ref{Thm_M_est_uniq} (with $X:=\cX_f=\RR$), we get that there exists $n_0\in\NN$ such that $\psi$ is not a $T_n$-function for any $n\geq n_0$, $n\in\NN$.
In particular, it yields that there exists $n_0\in\NN$ such that for each $n\geq n_0$, there exist real numbers $x_1,\ldots,x_n\in\RR$ such that the likelihood equation \eqref{psi_est_equation} based on $\bx=(x_1,\ldots,x_n)$ has at least two solutions.
In such a case, by part (v) of Theorem \ref{Thm_M_est_uniq}, we also get that there exists a  $\blambda\in\Lambda_2$ such that $\psi$ is not a $T_2^\blambda$-function.
\proofend
\end{Ex}

\section{Proofs for Section \ref{Section_Ex_uniq_empirical}}
\label{Section_Proof_for_Section2}

{\bf Proof of Lemma \ref{Lem_level_of_increase_and-sign_change}.}
Assume that $y$ is a level of increase for $f$ and (i) and (ii) are not valid. Define
 \[
  A:=\{v\in \Theta\mid y\geq f(v)\} \qquad\mbox{and}\qquad
  B:=\{u\in \Theta\mid y\leq f(u)\}.
 \]
Then $A\cup B=\Theta$ and $A,B$ are nonempty, since (i) and (ii) do not hold.
If $v\in A$ and $u\in B$, then $v\leq u$ (since $y$ is a level of increase for $f$).
Consequently, $\sup A\leq \inf B$, and, using that $A\cup B=\Theta$, we have that $\sup A=\inf B=:\vartheta$.
If $t<\vartheta$, then $t\not\in B$, which implies that $y>f(t)$, i.e., $y-f(t)>0$.
Similarly, if $\vartheta<t$, then $t\not\in A$, which yields $y<f(t)$, i.e., $y-f(t)<0$.
Hence, we have proved that $\vartheta$ is a point of sign change for the function $y-f$, i.e., (iii) must hold.

Conversely, if (i) or (ii) hold, then, just
using the definition, we have that $y$ is a level of increase for $f$.
Finally, assume that (iii) is valid, i.e., there exists $\vartheta\in \Theta$ which is a point of sign change for the function $y-f$, and, on the contrary, suppose that $y$ is not a level of increase for $f$.
Then there exist $u,v\in \Theta$ such that $u<v$ and $f(v)\leq y\leq f(u)$.
Therefore, $y-f(v)\geq0\geq y-f(u)$, which, using the definition of a point of sign change, implies
 that $v\leq\vartheta\leq u$ contradicting $u<v$.
\proofend

{\bf Proof of Lemma \ref{Lem_level_of_increase}.}
Assume that the levels of increase for $f$ form a dense subset in the convex hull of $f(\Theta)$, but $f$ is not increasing.
Then there exist $u,v\in \Theta$ with $u<v$ such that $f(u)>f(v)$.
The convex hull of $f(\Theta)$ contains the open interval $(f(v),f(u))$.
Hence, by the assumption, one can find an element $y\in (f(v),f(u))$ which is a level of increase for $f$. Therefore, by definition, $v\leq u$, which contradicts $u<v$.
The second statement of the lemma readily follows from the definitions of strictly increasing property and level of increase.
For the last statement of the lemma, assume that $u,v\in \Theta$ satisfy the inequalities
\[
g(f(v))=(g\circ f)(v)
\leq g(y)\leq (g\circ f)(u) = g(f(u)).
\]
Since $g$ is strictly increasing, we get that  $f(v)\leq y\leq f(u)$. Using that $y$ is a level of increase for $f$, it follows that $v\leq u$, implying the third statement of the lemma.
\proofend

{\bf Proof of Lemma \ref{Lem_e_monoton}.}
To verify the $\vare$-increasingness property of $f$,
let $u,v\in \Theta$ be arbitrary such that $u<v$.
The intervals $J_1:=[y_0,y_1]$, \dots, $J_n:=[y_{n-1},y_{n}]$ cover the image $f(\Theta)$.
Therefore, for some $i\in\{1,\dots,n\}$, we have that $f(u)\in J_i$, which implies $y_{i-1}\leq f(u)\leq y_i$.
Since $y_{i-1}$ is a level of increase for $f$, by part (iii) of Remark \ref{Rem_level_of_increase}, we have that $y_{i-1}<f(v)$.
On the other hand, due to the definition of $\vare$, it follows that $y_i\leq y_{i-1}+\vare$.
Thus, we get
 \[
 f(u)\leq y_i\leq y_{i-1}+\vare<f(v)+\vare,
 \]
 which was to be proved.
\proofend

{\bf Proof of Theorem \ref{Thm_M_est_uniq}.}
(i): Let $x,y\in X$ be such that $\vartheta_1(x) <\vartheta_1(y)$.
Since $\psi\in\Psi[T_1](X,\Theta)$, for all $t\in(\vartheta_1(x),\vartheta_1(y))$,
 we have $\psi(x,t)<0$ and $\psi(y,t)>0$.

To the contrary, assume that $\frac{\lambda_2}{\lambda_1}$ is not a level of increase for the function \eqref{function_newhanyados0}.
Then there exist $u,v\in\RR$ such that $\vartheta_1(x) < u< v < \vartheta_1(y)$ and
\[
      -\frac{\psi(x,v)}{\psi(y,v)} \leq  \frac{\lambda_2}{\lambda_1} \leq -\frac{\psi(x,u)}{\psi(y,u)}.
\]
Rearranging these inequalities, we get that
\[
  \lambda_1\psi(x,u)+\lambda_2\psi(y,u)\leq 0\leq \lambda_1\psi(x,v)+\lambda_2\psi(y,v).
\]
In view of the property $[T_2^{(\lambda_1,\lambda_2)}]$ of $\psi$, this implies that
 \[
  v\leq\vartheta_{2,\psi}^{(\lambda_1,\lambda_2)}(x,y)\leq u,
 \]
 which contradicts the inequality $u<v$.

To the contrary, assume that $\frac{\lambda_1}{\lambda_2}$ is not a level of increase for the function \eqref{function_newhanyados0}.
Similarly as before, we have that there exist $u,v\in\RR$ such that $\vartheta_1(x) < u< v < \vartheta_1(y)$ and
 \[
  \lambda_1\psi(y,u)+\lambda_2\psi(x,u)\leq 0\leq \lambda_1\psi(y,v)+\lambda_2\psi(x,v).
\]
In view of the property $[T_2^{(\lambda_1,\lambda_2)}]$ of $\psi$, this implies that
 $v\leq\vartheta_{2,\psi}^{(\lambda_1,\lambda_2)}(y,x)\leq u$, which again contradicts the inequality $u<v$.

(ii): Let $\psi\in\Psi[T_n^{(\lambda_1,\ldots,\lambda_n)}](X,\Theta)$ for some $n\in\NN\setminus\{1\}$ and $(\lambda_1,\ldots,\lambda_n)\in(0,\infty)^n$.
Then for each $k\in\{1,\ldots,n-1\}$, we have $\psi$ is a $T_2^{(\sum_{i=1}^k \lambda_i,\sum_{i=k+1}^n \lambda_i)}$-function,
 since for each $x,y\in X$ and $k\in\{1,\ldots,n-1\}$, it holds that
 \[
   \left( \sum_{i=1}^k \lambda_i\right)\psi(x,t) + \left( \sum_{i=k+1}^n \lambda_i\right)\psi(y,t)
      = \sum_{i=1}^n \lambda_i \psi(x_i,t),\qquad t\in \Theta,
 \]
where $x_i:=x$, $i\in\{1,\ldots,k\}$ and $x_i:=y$, $i\in\{k+1,\ldots,n\}$.
Consequently, the assertion readily follows from part (i) of the present theorem.

(iii): If $\psi$ is a $T_n$-function for some $n\in\NN\setminus\{1\}$, then it is a $T_n^{(\lambda_1,\ldots,\lambda_n)}$-function with $\lambda_1=\cdots=\lambda_n:=1$, thus the assertion follows from part (ii) of the present theorem.

(iv): Let $(n_i)_{i\in\NN}\subseteq\NN$ be a strictly increasing sequence such that
$\psi\in\Psi[T_{n_i}](X,\Theta)$ for all $i\in\NN$.
Let $x,y\in X$ such that $\vartheta_1(x) <\vartheta_1(y)$.
Then, by assertion (iii) of this theorem we have that the numbers
 \[
   \left\{\frac{k}{n_i-k}\Mid k\in\{1,\dots,n_i-1\}, i\in\NN \right\}
    = \left\{ \frac{1}{n_i-1}, \frac{2}{n_i-2}, \ldots, \frac{n_i-2}{2}, n_i-1 \Mid  i\in\NN \right\}
  \]
are levels of increase for the function \eqref{function_newhanyados0}.
We are going to apply Lemma \ref{Lem_level_of_increase}.
The convex hull of the range of the function \eqref{function_newhanyados0} is contained in $(0,\infty)$,
 so if we check that the set $\{\frac{k}{n_i-k}\mid k\in\{1,\dots,n_i-1\}, i\in\NN \}$ is dense in $(0,\infty)$,
 then Lemma \ref{Lem_level_of_increase} will imply that the function \eqref{function_newhanyados0} is increasing.
Since the function $g:(0,\infty)\to(0,1)$, $g(u):=\frac{u}{u+1}$, $u>0$, is bijective,
 it is enough to check that the set
 \[
  \left\{g\left(\frac{k}{n_i-k}\right)\Mid k\in\{1,\dots,n_i-1\}, i\in\NN \right\}
    = \left\{\frac{k}{n_i}\Mid k\in\{1,\dots,n_i-1\}, i\in\NN \right\}
 \]
 is dense in $g((0,\infty))=(0,1)$.
This readily follows, since if $(a,b)\subseteq (0,1)$ is an open interval, then
 there exists $i_0\in\NN$ such that $\frac{1}{n_{i_0}}<b-a$ (due to $n_i\to\infty$ as $i\to\infty$),
 and hence there exists $k\in\{1,\ldots,n_{i_0}-1\}$ such that $\frac{k}{n_{i_0}}\in(a,b)$.

Now we turn to prove the second statement of the assertion (iv).
Let $\ell,m\in\NN$.
We show that $\frac{\ell}{m}$ is a level of increase for the function  \eqref{function_newhanyados0}.
By assumption, there exists $i_0\in\NN$ such that $m+\ell$ divides $n_{i_0}$ and $\psi \in\Psi[T_{n_{i_0}}](X,\Theta)$.
By part (iii) of this theorem, we have that the elements of the set
 \[
   \left\{\frac{k}{n_{i_0}-k}\Mid k\in\{1,\dots,n_{i_0}-1\}\right\}
 \]
are levels of increase for the function \eqref{function_newhanyados0}.
By choosing $k:=\frac{\ell n_{i_0}}{m+\ell}\in\{1,\ldots,n_{i_0}-1\}$, we have
 \[
   \frac{k}{n_{i_0}-k} = \frac{\ell n_{i_0}}{ (m+\ell)n_{i_0} - \ell n_{i_0}}
                       = \frac{\ell}{m},
 \]
 yielding that $\frac{\ell}{m}$ is a level of increase for the function \eqref{function_newhanyados0}, as desired.

(v): Assume that $\psi\in\Psi[T_2^{\blambda}](X,\Theta)$ for each $\blambda\in\Lambda_2$.
Let $x,y\in X$ be such that $\vartheta_1(x)<\vartheta_1(y)$.
By assertion (i), it follows that $\lambda_1/\lambda_2$ is a level of increase
 for the function \eqref{function_newhanyados0} for each $\lambda_1,\lambda_2>0$.
Since $\lambda_1$ and $\lambda_2$ are arbitrary positive numbers, we get that each positive number is a level of increase for the positive function \eqref{function_newhanyados0}.
In view of the second statement of Lemma \ref{Lem_level_of_increase}, this implies that the function \eqref{function_newhanyados0}
is strictly  increasing.

(vi): This assertion is an immediate consequence of Theorem~\ref{Thm_theor_M_est_uniq} as stated in Corollary \ref{Cor_M_est_uniq+},
 therefore its proof is omitted here.
\proofend

{\bf Proof of Corollary \ref{Cor_M_est_uniq}.}
If (i) holds, then part (vi) of Theorem \ref{Thm_M_est_uniq} implies that (iii) is valid as well.
If (iii) holds, then (ii) is readily satisfied.
Finally, if (ii) holds, then part (v) of Theorem \ref{Thm_M_est_uniq} implies the validity of (i).
\proofend

{\bf Proof of Proposition \ref{Pro_M_est_uniq}.}
(i): First, let us suppose that for each $x\in X$, the function $\Theta\ni t\mapsto \psi(x,t)$ is decreasing.
Let $\vartheta_1(x)<s<t<\vartheta_1(y)$.
Then, since $\psi$ is a $T_1$-function, we have
 \[
  0>\psi(x,s)\geq \psi(x,t)
   \qquad \text{and}\qquad
  \psi(y,s)\geq \psi(y,t)>0.
 \]
Consequently, we get
 \[
   0< -\psi(x,s)\psi(y,t) \leq -\psi(x,t)\psi(y,s),
 \]
 which is equivalent to
 \[
  -\frac{\psi(x,s)}{\psi(y,s)}
      \leq -\frac{\psi(x,t)}{\psi(y,t)},
 \]
 yielding that the function  \eqref{function_newhanyados0}  is increasing.
The case when the function $\Theta\ni t\mapsto \psi(x,t)$ is strictly decreasing for each $x\in X$ can be handled similarly.

(ii): Let us suppose that for each $x\in X$, the function $\Theta\ni t\mapsto \psi(x,t)$ is strictly decreasing.
Let $n\in\NN$, $\blambda=(\lambda_1,\ldots,\lambda_n)\in\Lambda_n$ and $(x_1,\ldots,x_n)\in X^n$.
Since $\lambda_i\geq 0$, $i\in\{1,\ldots,n\}$ and $\lambda_1+\cdots +\lambda_n>0$, we get that
 the function $\Theta\ni t\mapsto \sum_{i=1}^n \lambda_i \psi(x_i,t)$ is strictly decreasing.
Using that $\psi\in\Psi[T_1](X,\Theta)$, we have that
 \[
   \sum_{i=1}^n \lambda_i \psi(x_i,t)
     \begin{cases}
          >0 & \text{if $t<\min(\vartheta_1(x_1),\ldots,\vartheta_1(x_n))$,}\\
          <0 & \text{if $t>\max(\vartheta_1(x_1),\ldots,\vartheta_1(x_n))$.}
    \end{cases}
 \]
Consequently, we have
 \begin{align*}
  &t^*:=\sup\Big\{ t\in \Theta :  \sum_{i=1}^n \lambda_i \psi(x_i,t)  >0 \Big\} \leq  \max(\vartheta_1(x_1),\ldots,\vartheta_1(x_n)),\\
  &t_*:=\inf\Big\{ t\in \Theta :  \sum_{i=1}^n \lambda_i \psi(x_i,t)  <0 \Big\} \geq  \min(\vartheta_1(x_1),\ldots,\vartheta_1(x_n)),
 \end{align*}
yielding that $t^*\leq t_*$.
Using the definition of infimum, supremum and that the map $\Theta\ni t\mapsto \sum_{i=1}^n \lambda_i \psi(x_i,t)$ is strictly decreasing,
 we get $t^*=t_*$.
Indeed, if $t^*<t_*$ were true, then, by the definition of infimum and supremum, $\sum_{i=1}^n \lambda_i \psi(x_i,t)=0$, $t\in(t^*,t_*)$,
 would hold, contradicting the strictly decreasing property of the function $\Theta\ni t\mapsto \sum_{i=1}^\lambda \lambda_i \psi(x_i,t)$.
All in all, we get $\vartheta_{n,\psi}^{\pmb{\lambda}}(x_1,\ldots,x_n) = t^*=t_*$, and then $\psi$ possesses the property $[T_n^{\pmb{\lambda}}]$.
\proofend

{\bf Proof of Proposition \ref{Pro_3}.}
Let $x,y\in X$ be such that $\vartheta_1(x)<\vartheta_1(y)$.
Since $\psi\in\Psi[T_1](X,\Theta)$, for all $t\in(\vartheta_1(x),\vartheta_1(y))$,
 we have $\psi(x,t)<0$ and $\psi(y,t)>0$, and hence the function \eqref{function_newhanyados} takes values in $(0,1)$.
By part (iii) of Theorem \ref{Thm_M_est_uniq}, we have that
 the elements of the set $\{\frac{k}{n-k}\mid k\in\{1,\dots,n-1\}\}$ are levels of increase for the function \eqref{function_newhanyados0}.
Note that
 \[
    \frac{\psi(x,t)}{\psi(x,t)-\psi(y,t)} = g\left(-\frac{\psi(x,t)}{\psi(y,t)}\right),
     \qquad t\in(\vartheta_1(x),\vartheta_1(y)),
 \]
 where $g:(0,\infty)\to\RR$, $g(u):=\frac{u}{u+1}$, $u>0$.
Since $g$ is strictly increasing and $(0,\infty)$ contains the range of the function \eqref{function_newhanyados0},
 by Lemma \ref{Lem_level_of_increase}, we get that the elements of the set
 \[
  \left\{g\left(\frac{k}{n-k}\right)\Mid k\in\{1,\dots,n-1\}\right\}
   = \left\{\frac{k}{n}\Mid k\in\{1,\dots,n-1\}\right\}
 \]
are levels of increase for the function \eqref{function_newhanyados}.
Since the function  \eqref{function_newhanyados0} is positive, we readily have that $0$ is also
 a level of increase for the function  \eqref{function_newhanyados0}.
Consequently, using that the function \eqref{function_newhanyados} takes values in $(0,1)$,
 the conditions of Lemma \ref{Lem_e_monoton} are satisfied with the choices $y_k:=\frac{k}{n}$, $k\in\{0,1,\ldots,n\}$,
 and then we get that the function \eqref{function_newhanyados} is strictly $\vare$-increasing with
 $\vare:=\max\Big\{\frac{k}{n} - \frac{k-1}{n}: k\in\{1,\ldots,n\}\Big\} = \frac{1}{n}$, as desired.
\proofend

{\bf Proof of Proposition \ref{T_n}.}
Now let us suppose that $m\in\{1,\ldots,n\}$ and $m$ is a divisor of $n$. Then there exists $k\in\NN$ such that $n=km$.
Consequently, for each $(y_1,\ldots,y_m)\in X^m$,
 with the notation
  \begin{align*}
    (x_1,\ldots,x_n):= (\underbrace{y_1,\ldots,y_1}_{k}, \underbrace{y_2,\ldots,y_2}_{k},\ldots, \underbrace{y_m,\ldots,y_m}_{k})\in X^n,
  \end{align*}
using that $\psi$ is a $T_n$-function, we have that
 \begin{align*}
  k\sum_{i=1}^m\psi(y_i,t)
   = \sum_{i=1}^n \psi(x_i,t)
    \begin{cases}
        >0 & \text{if $t<\vartheta_n(x_1,\ldots,x_n)$,}\\
        <0 & \text{if $t>\vartheta_n(x_1,\ldots,x_n)$.}
     \end{cases}
 \end{align*}
Hence $\psi\in\Psi[T_m](X,\Theta)$ with $\vartheta_m(y_1,\ldots,y_m):=\vartheta_n(x_1,\ldots,x_n)$.
\proofend

{\bf Proof of Proposition \ref{Pro_T_n_lambda}.}
For each $\pmb{y}:=(y_1,\ldots,y_m)\in X^m$ and $t\in \Theta$, we have
 \begin{align}\label{help_T_n_lambda}
   \psi_{\pmb{y},\pmb{\mu}}(t):=\sum_{\alpha=1}^m \mu_\alpha \psi(y_\alpha,t)
                               = \sum_{\alpha=1}^m \left( \sum_{i\in H_\alpha}\lambda_i \right) \psi(y_\alpha,t)
                               = \sum_{j=1}^n \lambda_j \psi(x_j,t)
                               = \psi_{\pmb{x},\pmb{\lambda}}(t),
 \end{align}
 where $\pmb{x}:=(x_1,\ldots,x_n)\in X^n$ is such that $x_j:=y_\alpha$ if $j\in H_\alpha$.
By the assumption, the value $\vartheta_n^{\pmb{\lambda}}(x_1,\ldots,x_n)$ is a point of sign change
for the function $\psi_{\pmb{x},\pmb{\lambda}}$.
Therefore, by \eqref{help_T_n_lambda}, we can see that
 the function $\psi_{\pmb{y},\pmb{\mu}}$ has the same point of sign change, and hence we have
 \[
  \vartheta_m^{\pmb{\mu}}(y_1,\ldots,y_m) = \vartheta_n^{\pmb{\lambda}}(x_1,\ldots,x_n),
 \]
yielding that $\psi\in\Psi[T^{\bmu}_m](X,\Theta)$.
\proofend

{\bf Proof of Proposition~\ref{Pro_BajType_est}.}
First, we check that $\psi\in\Psi[T_1](X,\Theta)$ with $\vartheta_1=f^{(-1)}\circ\varphi$.
Let $x\in X$ be fixed.
If $t<(f^{(-1)}\circ\varphi)(x)$, $t\in\Theta$, then $\psi(x,t)>0$, since otherwise
 $\psi(x,t)\leq 0$ would yield that $\varphi(x)\leq f(t)$,
 and hence, by Lemma \ref{Lem_gen_left_inverse}, we would have that
 $(f^{(-1)}\circ \varphi)(x)\leq (f^{(-1)}\circ f)(t)=t$, leading us to a contradiction.
Similarly, if $t>(f^{(-1)}\circ\varphi)(x)$, $t\in\Theta$, then $\psi(x,t)<0$, since otherwise $\psi(x,t)\geq 0$ would yield that $\varphi(x)\geq f(t)$, and hence, by Lemma \ref{Lem_gen_left_inverse}, we would have that
 $(f^{(-1)}\circ \varphi)(x)\geq (f^{(-1)}\circ f)(t)=t$, leading us to a contradiction.
All in all, for each $x\in X$, we have that
 \[
  \psi(x,t) \begin{cases}
              >0 & \text{if $t<(f^{(-1)}\circ\varphi)(x)$, $t\in\Theta$,}\\
              <0 & \text{if $t>(f^{(-1)}\circ\varphi)(x)$, $t\in\Theta$,}
            \end{cases}
 \]
 as desired.

Consequently, using also that for each $x\in X$, the function $\Theta\ni t\mapsto \psi(x,t)=p(x)(\varphi(x)-f(t))$ is strictly decreasing,
 part (ii) of Proposition \ref{Pro_M_est_uniq} implies that $\psi\in\Psi[T_n^{\pmb{\lambda}}](X,\Theta)$ for each $n\in\NN$ and $\pmb{\lambda}\in\Lambda_n$.

It remains to check that \eqref{help14} holds.
First, note that the right hand side of \eqref{help14} is well-defined, since
 \begin{align*}
   &\frac{\lambda_1p(x_1)\varphi(x_1)+\dots+\lambda_np(x_n)\varphi(x_n)}{\lambda_1p(x_1)+\dots+\lambda_np(x_n)} \\
   &\qquad = \frac{\lambda_1p(x_1)}{\lambda_1p(x_1)+\dots+\lambda_np(x_n)}\varphi(x_1)
         + \cdots + \frac{\lambda_np(x_n)}{\lambda_1p(x_1)+\dots+\lambda_np(x_n)} \varphi(x_n)\\
   &\qquad \in \conv(\varphi(X))\subseteq \conv(f(\Theta)),
 \end{align*}
 and $f^{(-1)}$ is defined on $\conv(f(\Theta))$ (see Lemma \ref{Lem_gen_left_inverse}).
Let $n\in\NN$, $(x_1,\ldots,x_n)\in X^n$ and $\blambda\in\Lambda_n$ be fixed.
If
 \[
   t<f^{(-1)}\bigg(\frac{\lambda_1p(x_1)\varphi(x_1)+\dots+\lambda_np(x_n)\varphi(x_n)}{\lambda_1p(x_1)+\dots+\lambda_np(x_n)}\bigg),
     \quad t\in\Theta,
 \]
 then $\sum_{i=1}^n \lambda_i\psi(x_i,t)>0$, since otherwise $\sum_{i=1}^n \lambda_i\psi(x_i,t) \leq 0$ would yield that
 \[
   \frac{\lambda_1p(x_1)\varphi(x_1)+\dots+\lambda_np(x_n)\varphi(x_n)}{\lambda_1p(x_1)+\dots+\lambda_np(x_n)}
      \leq f(t).
 \]
Hence, by Lemma \ref{Lem_gen_left_inverse}, we would have that
 \[
     f^{(-1)}\bigg(\frac{\lambda_1p(x_1)\varphi(x_1)+\dots+\lambda_np(x_n)\varphi(x_n)}{\lambda_1p(x_1)+\dots+\lambda_np(x_n)}\bigg)
           \leq (f^{(-1)}\circ f)(t)=t,
 \]
 leading us to a contradiction.
Similarly, we can easily see that the inequality
 \[
   t>f^{(-1)}\bigg(\frac{\lambda_1p(x_1)\varphi(x_1)+\dots+\lambda_np(x_n)\varphi(x_n)}{\lambda_1p(x_1)+\dots+\lambda_np(x_n)}\bigg),
     \quad t\in\Theta,
 \]
implies $\sum_{i=1}^n \lambda_i\psi(x_i,t)<0$.
These two properties together with that $\psi\in\Psi[T_n^\blambda](X,\Theta)$ yield the equality \eqref{help14}.
\proofend

\section{Proofs for Section \ref{Section_Ex_uniq_theoretical}}
\label{Section_Proof_for_Section3}

{\bf Proof of Theorem \ref{Thm_theor_M_est_uniq}.}
Define the sets $U,V\subseteq\Theta$  by
 \begin{align*}
  U:=\big\{s\in\Theta: \EE(\psi(\xi,s))\geq0\big\}\qquad  \text{and}\qquad V:=\big\{t\in\Theta: \EE(\psi(\xi,t))\leq 0\big\}.
 \end{align*}
Then, in view of assumption (v), we have that $s_0\in U$ and $t_0\in V$. In what follows, we show that $s\leq t$ holds for all $s\in U$ and $t\in V$.
To the contrary, assume that $t<s$, and, for any Borel subset $H\subseteq\Theta$, let us define
 \[
  \Omega_{H}:=\{\omega\in\Omega: \vartheta_1(\xi(\omega))\in H\}.
 \]
Then $\Omega_{H}\in \cA$ due to the measurability of $ \vartheta_1:X\to\Theta$ and $\xi:\Omega\to X$.
Indeed, for each $r\in\Theta$ we have that
 \[
  \vartheta_1^{-1}((-\infty,r)) = \big\{ x\in X : \vartheta_1(x)< r\big\}
                                = \big\{ x\in X : \psi(x,r) < 0 \big\}
                                \in \cX,
 \]
 where we used assumptions (i), (iii) and that the sigma-algebra generated by the family
$\{(-\infty,r)\cap\Theta, r\in\Theta\}$ coincides with the Borel sigma-algebra on $\Theta$.

Consider the following partition of $\Theta$, which is induced by $t$ and $s$:
\[
 I:=\Theta\cap(-\infty,t),\qquad
 J:=[t,s],\qquad
 K:=\Theta\cap(s,\infty).
\]
Then, using assumption (i), we have
\begin{align*}
 \Omega_I&=\{\omega\in\Omega:\vartheta_1(\xi(\omega))<t\}
 =\{\omega\in\Omega:\psi(\xi(\omega),t) < 0\},\\
 \Omega_K&=\{\omega\in\Omega:\vartheta_1(\xi(\omega))>s\}
 =\{\omega\in\Omega:\psi(\xi(\omega),s) > 0\}.
\end{align*}
We show that $\PP(\Omega_I)>0$ and $\PP(\Omega_K)>0$.
Indeed, on the contrary, if $\PP(\Omega_I)=0$, then $\PP(\psi(\xi,t)\geq0)=1$,
 which implies that $\EE(\psi(\xi,t))\geq0$. By the inclusion $t\in V$,
 we also have that $\EE(\psi(\xi,t))\leq0$ and hence $\EE(\psi(\xi,t))=0$.
Therefore, $\PP(\psi(\xi,t)=0)=1$, i.e., $\PP(\vartheta_1(\xi)=t)=1$.
It follows from the inequality $t<s$ that $\PP(\vartheta_1(\xi)<s)=1$, and hence $\PP(\psi(\xi,s)<0)=1$.
This implies that $\EE(\psi(\xi,s))<0$, which contradicts that $s$ belongs to $U$.
The equality $\PP(\Omega_K)=0$ leads to a contradiction similarly.

The inequalities $\PP(\Omega_I)>0$ and $\PP(\Omega_K)>0$ imply that $\Omega_I\ne \emptyset$ and $\Omega_K\ne\emptyset$. Then, for all $\omega'\in \Omega_I$ and $\omega''\in \Omega_K$, we have that
\[
   \vartheta_1(\xi(\omega'))<t<s<\vartheta_1(\xi(\omega'')).
\]
Therefore, using assumption (ii) with $x:=\xi(\omega')$ and $y:=\xi(\omega'')$,
 the function \eqref{function_newhanyados0} is strictly increasing, hence we get
 \begin{align*}
  \frac{\psi(\xi(\omega'),s)}{\psi(\xi(\omega''),s)}
  <\frac{\psi(\xi(\omega'),t)}{\psi(\xi(\omega''),t)}.
 \end{align*}
Using that $\psi(\xi(\omega''),s)>0$ and $\psi(\xi(\omega''),t)>0$, we can obtain that
 \begin{align}\label{help11}
  \psi(\xi(\omega'),s)\psi(\xi(\omega''),t)
  <\psi(\xi(\omega'),t)\psi(\xi(\omega''),s),
  \qquad(\omega',\omega'')\in \Omega_I\times\Omega_K.
 \end{align}
Integrating on $\Omega_I$ and then on $\Omega_K$ with respect to $\PP$, it follows that
 \[
  \int_{\Omega_I}\psi(\xi(\omega'),s)\,\dd\PP(\omega')\cdot\int_{\Omega_K}\psi(\xi(\omega''),t)\,\dd\PP(\omega'')
  <\int_{\Omega_I}\psi(\xi(\omega'),t)\,\dd\PP(\omega')\cdot\int_{\Omega_K}\psi(\xi(\omega''),s)\,\dd\PP(\omega''),
 \]
 that is,
 \begin{align}\label{help9}
  \EE(\psi(\xi,s)\bone_{\Omega_I})\cdot\EE(\psi(\xi,t)\bone_{\Omega_K})
  <\EE(\psi(\xi,t)\bone_{\Omega_I})\cdot\EE(\psi(\xi,s)\bone_{\Omega_K}).
 \end{align}
The inequality in \eqref{help9} is indeed strict because the left hand side of \eqref{help11} is strictly smaller than
 its right hand side over the set $\Omega_I\times\Omega_K$ which has positive measure with respect to the product probability $\PP\otimes\PP$.

Furthermore, using also that $t<s$, for each $\omega'\in \Omega_J$ and $\omega''\in\Omega_K$, we have that $\psi(\xi(\omega'),s)\leq0$, $\psi(\xi(\omega''),t)>0$, $\psi(\xi(\omega'),t)\geq0$, and $\psi(\xi(\omega''),s)>0$. Therefore,
\[
  \psi(\xi(\omega'),s)\psi(\xi(\omega''),t)
  \leq 0\leq\psi(\xi(\omega'),t)\psi(\xi(\omega''),s),
  \qquad(\omega',\omega'')\in \Omega_J\times\Omega_K.
\]
Integrating on $\Omega_J$ and then on $\Omega_K$ with respect to $\PP$, it follows that
 \begin{align}\label{help12}
  \EE(\psi(\xi,s)\bone_{\Omega_J})\cdot\EE(\psi(\xi,t)\bone_{\Omega_K})
      \leq\EE(\psi(\xi,t)\bone_{\Omega_J})\cdot\EE(\psi(\xi,s)\bone_{\Omega_K}).
 \end{align}
Adding up the inequalities \eqref{help9} and \eqref{help12}, and using that $\Omega_I$ and $\Omega_J$ are disjoint, we get
\[
  A(s)B(t):=\EE(\psi(\xi,s)\bone_{\Omega_I\cup\Omega_J})
             \cdot\EE(\psi(\xi,t)\bone_{\Omega_K})
           <\EE(\psi(\xi,t)\bone_{\Omega_I\cup\Omega_J} )\cdot\EE(\psi(\xi,s)\bone_{\Omega_K})=:A(t)B(s).
\]
Further, we have $A(t)+B(t) = \EE(\psi(\xi,t)) \leq 0$ and $A(s)+B(s) = \EE(\psi(\xi,s))\geq 0$, since $t\in V$, $s\in U$, and $\Omega_I\cup\Omega_J$ and $\Omega_K$ are disjoint.

To summarize, $t,s\in\Theta$ are such that $t<s$ and the following inequalities hold
\begin{align}\label{help13}
  A(s) B(t)<B(s)A(t), \qquad A(t)+B(t)\leq 0,
  \qquad A(s)+B(s)\geq 0.
\end{align}
Here $B(t)>0$ because it equals the integral of a positive function over the set $\Omega_K$ which has positive measure with respect to the probability $\PP$.
On the other hand, $A(s)<0$, because
 \[
  A(s)=\EE(\psi(\xi,s)\bone_{\Omega_I\cup\Omega_J})
  =\EE(\psi(\xi,s)\bone_{\Omega_I})+\EE(\psi(\xi,s)\bone_{\Omega_J})
 \]
 and the first term is negative being equal to the
 integral of a negative function over the set $\Omega_I$ (which has positive measure with respect to $\PP$)
 and the second term is nonpositive being equal to the integral of a nonpositive function over the set $\Omega_J$.

Consequently, by the last two inequalities of \eqref{help13}, we get
 \[
  0<B(t)\leq -A(t) \qquad\text{and}\qquad 0<-A(s)\leq B(s),
 \]
 yielding that
 \[
    0 < -A(s)B(t) \leq -A(t)B(s),
 \]
 i.e., $A(s)B(t)\geq B(s)A(t)$.
This contradicts to the first inequality in \eqref{help13}.

Consequently, we have that $s_0\leq u_0:=\sup U \leq \inf V=: v_0\leq t_0$. It remains to show that $u_0 = v_0$.
If, to the contrary, we assume that $u_0 < v_0$, then for each $r\in(u_0,v_0)$, we get $r\notin U$ and $r\notin V$,
 yielding that $\EE(\psi(\xi,r))< 0$ and $\EE(\psi(\xi,r))> 0$, respectively, which is a contradiction.

All in all, $u_0=v_0$ is a unique point of sign change for the function $\Theta\ni t\mapsto \EE(\psi(\xi,t))$,   as desired.
\proofend

{\bf Proof of Proposition \ref{Pro_theor_M_est_uniq2}.}
By the assumption (i), for each $\omega\in\Omega$, we have that the function $\Theta\ni t\mapsto \psi(\xi(\omega),t)$
 is strictly decreasing.
By the monotonicity of the expectation,  it implies that the function $\Theta\ni t\mapsto \EE(\psi(\xi,t))$
 is decreasing, and in fact, it is strictly decreasing.
Indeed, if $t_1<t_2$, $t_1,t_2\in\Theta$, are such that $\EE(\psi(\xi,t_1)) =  \EE(\psi(\xi,t_2))$, then
 $\EE(\psi(\xi,t_1) - \psi(\xi,t_2))=0$, where $\PP(\psi(\xi,t_1) - \psi(\xi,t_2)\geq 0)=1$.
Consequently, $\PP(\psi(\xi,t_1) - \psi(\xi,t_2) = 0)=1$, leading us to a contradiction, since
 $\psi(\xi(\omega),t_1)> \psi(\xi(\omega),t_2)$, $\omega\in\Omega$.

Define the sets $U,V\subseteq\Theta$  by
 \begin{align*}
  U:=\big\{s\in\Theta: \EE(\psi(\xi,s))\geq0\big\}\qquad  \text{and}\qquad V:=\big\{t\in\Theta: \EE(\psi(\xi,t))\leq 0\big\}.
 \end{align*}
By the assumption (iv), we have that $s_0\in U$ and $t_0\in V$.
Since the function $\Theta\ni t\mapsto \EE(\psi(\xi,t))$ is strictly decreasing, we can easily deduce that $s\leq t$ for all $s\in U$, $t\in V$.
Indeed, if for some $s\in U$ and $t\in V$, the inequality $s>t$ were true, then we would have that
 $0\leq \EE(\psi(\xi,s)) < \EE(\psi(\xi,t)) \leq 0$, leading us to a contradiction.
Hence we have $s_0\leq u_0:=\sup U\leq \inf V =: v_0\leq t_0$.
It remains to show that $u_0 = v_0$.
If, to the contrary, we assume that $u_0 < v_0$, then for each $r\in(u_0,v_0)$, we get $r\notin U$ and $r\notin V$,
 yielding that $\EE(\psi(\xi,r))< 0$ and $\EE(\psi(\xi,r))> 0$, respectively, which is a contradiction.
Consequently, $u_0=v_0$ is a unique point of sign change for the function $\Theta\ni t\to\EE(\psi(\xi,t))$, as desired.
\proofend

{\bf Proof of Corollary \ref{Cor_M_est_uniq+}.}
To verify the statement, we have to show that, for each $n\in\NN$,
 $x_1,\dots,x_n\in X$ and $\blambda=(\lambda_1,\dots,\lambda_n)\in\Lambda_n$,
 the function $\Theta\ni t\mapsto \sum_{i=1}^n\lambda_i\psi(x_i,t)$ has a unique point of sign change in $\Theta$.
Without loss of generality, we may assume that $x_1,\dots,x_n$ are  pairwise distinct elements of $X$ and
 $\lambda_1,\dots,\lambda_n>0$ with $\lambda_1+\dots+\lambda_n=1$.

Define the probability space $(\Omega,\cA,\PP)$ by
 \[
  \Omega:=\{x_1,\dots,x_n\}, \qquad \cA:=2^\Omega,\qquad
  \PP(\{x_i\}):=\lambda_i \qquad i\in\{1,\dots,n\},
 \]
 and the random variable $\xi:\Omega\to\Omega$ by $\xi(\omega):=\omega$, $\omega\in\Omega$.

Then the conditions (i) and (ii) of Theorem~\ref{Thm_theor_M_est_uniq} follow from our assumptions.
The measurability condition (iii) of Theorem~\ref{Thm_theor_M_est_uniq} is trivial due to the fact that $\cA=2^\Omega$.
Since
 \[
  \EE(|\psi(\xi,t)|)=\sum_{i=1}^n\lambda_i|\psi(x_i,t)|,\qquad t\in\Theta,
 \]
 the condition (iv) of Theorem~\ref{Thm_theor_M_est_uniq} is obviously valid. Finally, the condition (v) of Theorem~\ref{Thm_theor_M_est_uniq} is satisfied by
\[
 s_0:=\min\{\vartheta_1(x_1),\dots,\vartheta_1(x_n)\}\qquad\text{and}\qquad
 t_0:=\max\{\vartheta_1(x_1),\dots,\vartheta_1(x_n)\}.
\]
Indeed, for each $i\in\{1,\dots,n\}$, we have that $\psi(x_i,s_0)\geq0\geq\psi(x_i,t_0)$, since $\psi$ is a $T_1$-function and $\psi(x,\vartheta_1(x))=0$, $x\in X$.
This implies that
\[
   \EE(\psi(\xi,s_0))=\sum_{i=1}^n\lambda_i\psi(x_i,s_0)
   \geq0
   \geq \sum_{i=1}^n\lambda_i\psi(x_i,t_0)=\EE(\psi(\xi,t_0)).
\]
Therefore, according to the conclusion of Theorem~\ref{Thm_theor_M_est_uniq},
 the mapping $\Theta\ni t\mapsto \EE(\psi(\xi,t))=\sum_{i=1}^n\lambda_i\psi(x_i,t)$ has a unique point of sign change in $\Theta$,
 as desired.
\proofend

{\bf Proof of Lemma \ref{Lem_conv_hull}.}
Let us define the probability measure $\QQ$ on the measurable space $(\Omega,\cA)$ by
 \[
  \QQ(A):=\int_A \frac{p(\xi)}{\EE(p(\xi))}\,\dd \PP = \frac{\EE(p(\xi)\bone_{A})}{\EE(p(\xi))},\qquad A\in\cA.
 \]
By denoting the expectation with respect to $\QQ$ by $\EE_\QQ$, we have
 \[
    \EE_\QQ(\vert \varphi(\xi)\vert ) = \frac{\EE(p(\xi)\vert \varphi(\xi)\vert)}{\EE(p(\xi))},
 \]
 and hence, by the assumptions, $\varphi(\xi)$ is integrable with respect to $\QQ$, and we also get
 \[
    \EE_\QQ(\varphi(\xi)) = \frac{\EE(p(\xi)\varphi(\xi))}{\EE(p(\xi))}.
 \]
Applying Lemma 1 in Jankovi\'c and Merkle \cite{JanMer} (for integrable one-dimensional random variables),
 we have that $\EE_\QQ(\varphi(\xi))\in \conv(\varphi(\xi(\Omega)))\subseteq \conv(\varphi(X))$, yielding the statement.
\proofend

{\bf Proof of Proposition~\ref{Pro_Bajrak_theor}.}
We apply Proposition \ref{Pro_theor_M_est_uniq2}.
The assumptions (i), (ii) and (iii) of Proposition \ref{Pro_theor_M_est_uniq2} readily hold.

To verify the assumption (iv) of Proposition \ref{Pro_theor_M_est_uniq2},
 we first show that, for any  $y\in J:=\conv(f(\Theta))$,
 there exist $s_0,t_0\in\Theta$ such that $f(s_0)\leq y\leq f(t_0)$.
By the Carath\'eodory's Theorem on convex hulls,
 there exist at most two elements $y_1,y_2\in f(\Theta)\subseteq J$ with $y_1\leq y_2$ such that $y$ can be represented as a convex combination of $y_1$ and $y_2$.
This also yields that $y_1\leq y\leq y_2$, and therefore, there exist $s_0,t_0\in\Theta$ such that $f(s_0)\leq y\leq f(t_0)$.
Now observe that
 \[
  f^{(-1)}\bigg(\frac{\EE(p(\xi)\varphi(\xi))}{\EE(p(\xi))}\bigg)
 \]
is well-defined, since, by Lemma \ref{Lem_conv_hull}, we get that
 \[
   \frac{\EE(p(\xi)\varphi(\xi))}{\EE(p(\xi))}\in \conv(\varphi(X))\subseteq \conv(f(\Theta)),
 \]
and, by Lemma \ref{Lem_gen_left_inverse}, $f^{(-1)}$ is defined on $\conv(f(\Theta))$.
Next, for $y:=\EE(p(\xi)\varphi(\xi))/\EE(p(\xi))$, let us choose $s_0,t_0\in\Theta$ as it was described above. Then
\[
 \EE(\psi(\xi,s_0))
 =\EE(p(\xi)\varphi(\xi))-f(s_0)\EE(p(\xi))
 =\EE(p(\xi))(y-f(s_0))\geq0,
\]
and
\[
 \EE(\psi(\xi,t_0))
 =\EE(p(\xi)\varphi(\xi))-f(t_0)\EE(p(\xi))
 =\EE(p(\xi))(y-f(t_0))\leq0.
\]
Therefore, the assumption (iv) of Proposition \ref{Pro_theor_M_est_uniq2} holds as well, and,
according to the conclusion of Proposition \ref{Pro_theor_M_est_uniq2},
we get that the function $\Theta\ni t\to\EE(\psi(\xi,t))$ admits a unique point of sign change in $\Theta$.
It remains to check that this unique point of sign change takes the form given in the proposition.

If, for some $t\in\Theta$, we have $t<f^{(-1)}(y)$, then $\EE(\psi(\xi,t))>0$, since otherwise $\EE(\psi(\xi,t))\leq 0$ would yield that
 $\EE(p(\xi)\varphi(\xi))\leq f(t)\EE(p(\xi))$, i.e., $y\leq f(t)$.
Then, by Lemma \ref{Lem_gen_left_inverse}, we would get
$f^{(-1)}(y) \leq f^{(-1)}(f(t)) = t$ leading us to a contradiction.

If for some $t\in\Theta$, we have $t>f^{(-1)}(y)$,
then one can similarly argue to obtain that $\EE(\psi(\xi,t))<0$.

Consequently, the unique point of sign change in question is $f^{(-1)}(y)$, as desired.
\proofend

{\bf Proof of Proposition \ref{Pro_Math_type}.}
First, we give a direct proof.
Denote the limit $\lim_{z\to\infty} f(z)$ by $f_\infty \in(0,\infty)$.
In view of the increasingness of $f$, it follows that $0\leq f(z)\leq f_\infty$ for all $z\in\RR_+$. Therefore,
$|\widetilde{f}(z)|\leq f_\infty$ for all $z\in\RR$, which implies that $|\psi(x,t)|\leq f_\infty$ for all $x,t\in\RR$.
Hence, for any random variable $\xi$ and for any $t\in\RR$, we have that $\EE(|\psi(\xi,t)|)<\infty$.

Since for each $x\in\RR$, the function \ $\RR\ni t\mapsto \psi(x,t)$ \ is strictly decreasing, we have
 the function $\RR\ni t\mapsto \EE(\psi(\xi,t))$ is strictly decreasing.
Indeed, if $s<t$, $s,t\in\RR$, then we have $\psi(\xi(\omega),s)>\psi(\xi(\omega),t)$, $\omega\in\Omega$,
 yielding that $\EE(\psi(\xi,s))\geq \EE(\psi(\xi,t))$.
Here the equality cannot hold, since otherwise $\EE(\psi(\xi,s) - \psi(\xi,t))=0$ would be valid
 yielding that $\PP(\psi(\xi,s)-\psi(\xi,t)=0)=1$.
This leads us to a contradiction, since $\psi(\xi(\omega),s)-\psi(\xi(\omega),t)>0$, $\omega\in\Omega$.
Using that
 $\lim_{t\to\pm\infty}\psi(\xi(\omega),t) = \mp f_\infty$, $\omega\in\Omega$,
 and $\vert \psi(\xi(\omega),t)\vert\leq f_\infty$, $\omega\in\Omega$, $t\in\RR$,
 the dominated convergence theorem implies that
 \begin{align}\label{help_Mat_1}
   \lim_{t\to\infty} \EE(\psi(\xi,t))
     = \EE\Big(\lim_{t\to\infty} \psi(\xi,t)\Big)
     = \EE(-f_{\infty})= -f_{\infty}<0,
 \end{align}
 and
 \begin{align}\label{help_Mat_2}
  \lim_{t\to-\infty} \EE(\psi(\xi,t))
   = \EE\Big(\lim_{t\to-\infty} \psi(\xi,t)\Big)
   = \EE(f_{\infty})= f_{\infty}>0.
 \end{align}
Since $f$ is continuous and $f(0)=0$, we have that  $\psi$ is continuous in its second variable.
Thus, by the dominated convergence theorem, it follows that the function $\RR\ni t\mapsto \EE(\psi(\xi,t))$ is also continuous.
All in all, the function $\RR\ni t \mapsto \EE(\psi(\xi,t))$ is strictly decreasing, continuous, and changes sign,
 and hence there exists a unique $t_0\in\RR$ such that $\EE(\psi(\xi,t_0))=0$, as desired.

Finally, we present an alternative proof of Proposition \ref{Pro_Math_type} using Theorem \ref{Thm_theor_M_est_uniq}.
We check that the assumptions of Theorem \ref{Thm_theor_M_est_uniq} hold.
Since $f(0)=0$ and $f$ is strictly increasing, the assumption (i) of Theorem \ref{Thm_theor_M_est_uniq} holds with $\vartheta_1(x)=x$, $x\in\RR$.
Using that $f$ is strictly increasing, by part (d) of  Proposition \ref{Pro_Mat}, we have that $\psi$ is a $T_2^{\blambda}$-function for all $\blambda\in\Lambda_2$.
Consequently, part (v) of Theorem \ref{Thm_M_est_uniq} yields that, for each $x,y\in\RR$ with $\vartheta_1(x)<\vartheta_1(y)$,
 the function \eqref{function_newhanyados0} is strictly increasing,
 i.e., the assumption (ii) of Theorem \ref{Thm_theor_M_est_uniq} holds.
The assumption (iii) of Theorem \ref{Thm_theor_M_est_uniq} readily holds.
The first part of the direct proof of the present proposition implies that the assumption (iv) of Theorem
 \ref{Thm_theor_M_est_uniq} holds.
Using \eqref{help_Mat_1} and \eqref{help_Mat_2} we have that the assumption (v) of Theorem \ref{Thm_theor_M_est_uniq} holds as well.
All in all, we can apply Theorem \ref{Thm_theor_M_est_uniq}, and it yields that
 the function $\RR\ni t\mapsto \EE(\psi(\xi,t))=\EE(\sign(\xi-t)f(\vert \xi-t\vert))$ has a (unique) point of sign change.
Since the function $\RR\ni t\mapsto \EE(\psi(\xi,t))$ is strictly decreasing and continuous (see the direct proof),
 we have that the equation $\EE(\psi(\xi,t))=0$ has a unique solution with respect to $t\in\RR$, as desired.
 \proofend

\section{Proofs for Section \ref{Section_statistical_examples}}
\label{Section_Proof_for_Section4}

{\bf Proof of Proposition \ref{Pro_quantile}.}
Assume that $n\geq 2$ and $\psi\in\Psi[T_n](\RR,\RR)$.
Then, according to part (iii) of Theorem \ref{Thm_M_est_uniq}, for each $x,y\in\RR$ with $x<y$,
 the number $\frac{n-k}{k}$ must be a level of increase
 for the function \eqref{function_newhanyados0} if $k\in\{1,\dots,n-1\}$.
Note that for each $x,y\in\RR$ with $x<y$, the function \eqref{function_newhanyados0} takes the form
 \[
   (x,y)\ni t\mapsto -\frac{\psi(x,t)}{\psi(y,t)} = \frac{1-\alpha}{\alpha} >0.
 \]
Therefore, we have that
 \[
 \frac{1-\alpha}{\alpha}\neq \frac{n-k}{k} \qquad \text{for each $k\in\{1,\dots,n-1\}$,}
 \]
 which implies that $\alpha\neq\frac{k}{n}$ for each $k\in\{1,\dots,n-1\}$.

Conversely, if $\alpha\not\in\big\{\frac{1}{n},\dots,\frac{n-1}{n}\big\}$, where $n\geq 2$,
 then we have that
 \begin{align}\label{help20}
  \frac{k}{n}<\alpha<\frac{k+1}{n}
 \end{align}
for some $k\in\{0,\dots,n-1\}$, and hence $k<n\alpha<k+1$.
Let $x_1,\dots,x_n\in \RR$ be arbitrary. If $t\in\RR$ with $t<x^*_{k+1}$, then we have that
$\psi(x^*_i,t)\geq \alpha-1$, $i=1,\ldots,k$, and $\psi(x^*_i,t)=\alpha$, $i=k+1,\ldots,n$, yielding that
\[
\sum_{i=1}^n \psi(x_i,t)
=\sum_{i=1}^k \psi(x^*_i,t) +\sum_{i=k+1}^n \psi(x^*_i,t)
  \geq k(\alpha-1)+(n-k)\alpha=n\alpha-k>0.
\]
If $t\in\RR$ with $t>x^*_{k+1}$, then we have that $\psi(x^*_i,t) = \alpha-1$, $i=1,\ldots,k+1$, and $\psi(x^*_i,t)\leq \alpha$, $i=k+2,\ldots,n$, yielding that
\[
\sum_{i=1}^n \psi(x_i,t)
=\sum_{i=1}^{k+1} \psi(x^*_i,t) +\sum_{i=k+2}^n \psi(x^*_i,t)
  \leq (k+1)(\alpha-1)+(n-k-1)\alpha=n\alpha-k-1<0.
\]
Therefore $\vartheta_n(x_1,\dots,x_n)$ exists and equals $x^*_{k+1}$.
This proves that $\psi$ is indeed a $T_n$-function.
Furthermore, using \eqref{help20}, we have $k<\alpha n<k+1$ and $k+1<\alpha n +1< k+2$, yielding that
 $\lceil n\alpha \rceil = \lfloor n\alpha+1 \rfloor=k+1$.
Hence $x_{k+1}^*=x^*_{\lceil n\alpha \rceil}=x^*_{\lfloor n\alpha+1 \rfloor}$, as desired.
\proofend

{\bf Proof of Proposition \ref{Pro_expectile}.}
Let us apply Proposition \ref{Pro_theor_M_est_uniq2} with the choices $X:=\RR$, $\cX:=\cB(\RR)$ and $\Theta:=\RR$.
The assumptions (i) and (ii) of Proposition \ref{Pro_theor_M_est_uniq2} readily hold.

The validity of the assumption (iii) of Proposition \ref{Pro_theor_M_est_uniq2} can be seen from
 \begin{align}\label{help_expectile2}
  \begin{split}
  \EE(\vert \psi(\xi,t)\vert)
    &= \EE(\alpha \vert \xi -t \vert \bone_{\{\xi>t\}} +  (1-\alpha) \vert \xi -t \vert \bone_{\{\xi<t\}} )
     \leq \alpha \EE(\vert \xi -t \vert) + (1-\alpha)\EE(\vert \xi -t \vert)\\
    &\leq \EE(\vert \xi \vert) + \vert t \vert<\infty,\qquad t\in\RR.
 \end{split}
 \end{align}

Finally, we verify the assumption (iv) of Proposition \ref{Pro_theor_M_est_uniq2}.
First, note that for all $t\in\RR$, we have
 \begin{align}\label{help_expectile-1}
  \begin{split}
  \EE(\psi(\xi,t))
   & = \alpha \EE( (\xi-t)^+ )  - (1-\alpha) \EE( (\xi-t)^{-})\\
   & = \alpha \EE(\xi-t) + (2\alpha-1)\EE( (\xi-t)^{-} ),
 \end{split}
 \end{align}
and, analogously,
 \begin{align}\label{help_expectile-2}
  \begin{split}
  \EE(\psi(\xi,t))
   = (2\alpha -1) \EE( (\xi-t)^+ ) + (1-\alpha)\EE(\xi-t).
  \end{split}
 \end{align}

In case of $\alpha=\frac{1}{2}$, we have $\EE(\psi(\xi,t)) = \frac{1}{2}(\EE(\xi)-t)$,
 which is positive if $t<\EE(\xi)$, and is negative if $t>\EE(\xi)$.
This shows that the assumption (iv) of Proposition \ref{Pro_theor_M_est_uniq2} holds in case of $\alpha=\frac{1}{2}$.

In case of $\alpha\in(\frac{1}{2},1)$, we have $2\alpha-1>0$, and hence \eqref{help_expectile-1} yields that
 $\EE(\psi(\xi,t)) \geq \alpha (\EE(\xi) -t)$, which is positive if $t<\EE(\xi)$.
Further, for all $t>0$, we have
 \begin{align*}
  \EE(\psi(\xi,t))
  \leq \alpha(\EE(\xi) - t) + (2\alpha-1)(\EE(\vert \xi\vert) + t)
   = \alpha\EE(\xi) + (2\alpha-1)\EE(\vert \xi\vert) - (1-\alpha)t,
 \end{align*}
 which is negative if
 \[
   t > \max\left(0,\frac{\alpha\EE(\xi) + (2\alpha-1)\EE(\vert \xi\vert)}{1-\alpha}\right).
 \]
This shows that the assumption (iv) of Proposition \ref{Pro_theor_M_est_uniq2} holds in case of $\alpha\in(\frac{1}{2},1)$.

In case of $\alpha\in(0,\frac{1}{2})$, we have $2\alpha-1<0$, and hence \eqref{help_expectile-2} yields that
 $\EE(\psi(\xi,t)) \leq (1-\alpha) (\EE(\xi) -t)$, which is negative if $t>\EE(\xi)$.
Further, for all $t<0$, we have
 \begin{align*}
  \EE(\psi(\xi,t))
  \geq (2\alpha-1)(\EE(\vert \xi\vert) + \vert t\vert) + (1-\alpha)(\EE(\xi) - t)
   = (2\alpha-1)\EE(\vert \xi\vert) + (1-\alpha)\EE(\xi) - \alpha t,
 \end{align*}
 which is positive if
 \[
   t < \min\left(0,\frac{(2\alpha-1)\EE(\vert \xi\vert) + (1-\alpha)\EE(\xi)}{\alpha}\right).
 \]
This shows that the assumption (iv) of Proposition \ref{Pro_theor_M_est_uniq2} holds in case of $\alpha\in(0,\frac{1}{2})$.

Therefore, the assumption (iv) of Proposition \ref{Pro_theor_M_est_uniq2} holds as well, and,
according to the conclusion of Proposition \ref{Pro_theor_M_est_uniq2},
we get that the function $\RR\ni t\to\EE(\psi(\xi,t))$ admits a unique point of sign change.

Using the dominated convergence theorem, we check that the function
 $\RR\ni t\mapsto \EE(\psi(\xi,t))$ is continuous.
Let $(t_n)_{n\in\NN}$ be a real sequence such that $t_n\to t_0$ as $n\to\infty$, where $t_0\in\RR$.
Then, using that $\psi$ is strictly decreasing in its second variable, we have that
 \[
   \psi\Big(\xi,\sup_{m\in\NN}t_m\Big) \leq \psi(\xi,t_n) \leq \psi\Big(\xi,\inf_{m\in\NN}t_m\Big), \qquad n\in\NN,
 \]
 yielding that
 \[
   \EE\Big(\psi\Big(\xi,\sup_{m\in\NN}t_m\Big)\Big) \leq \EE(\psi(\xi,t_n)) \leq \EE\Big(\psi\Big(\xi,\inf_{m\in\NN}t_m\Big)\Big), \qquad n\in\NN.
 \]
Hence, by \eqref{help_expectile2}, we get that
 \[
  \EE(\vert \psi(\xi,t_n)\vert)\leq \EE\Big(\Big \vert \psi\Big(\xi,\inf_{m\in\NN}t_m\Big) \Big\vert \Big)
                                     + \EE\Big( \Big\vert \psi\Big(\xi,\sup_{m\in\NN}t_m\Big) \Big\vert \Big)
     <\infty, \qquad n\in\NN.
 \]
Further, using that  $\psi$ is continuous in its second variable, we have $\psi(\xi,t_n)\to\psi(\xi,t_0)$ as $n\to\infty$.
Hence the dominated convergence theorem implies that $\EE(\psi(\xi,t_n))\to \EE(\psi(\xi,t_0))$ as $n\to\infty$,
 as desired.

Consequently, the unique point of sign change of the function $\RR\ni t\to\EE(\psi(\xi,t))$ is nothing else but the unique solution of the equation $\EE(\psi(\xi,t))=0$, $t\in\RR$.
\proofend

{\bf Proof of Proposition \ref{Pro_Mat}.}
Part (a):
It follows from the facts that for each $x\in\RR$, we have $\psi(x,x)=0$; $\psi(x,t)>0$ holds for each $t<x$ if and only if $f(z)>0$ for each $z>0$;
 and  $\psi(x,t)<0$ holds for each $t>x$ if and only if $f(z)>0$ for each $z>0$.

Part (b):
Let us suppose that $\psi\in\Psi[T_n](\RR,\RR)$ for infinitely many $n\in\NN$.
Then, by Proposition \ref{T_n}, we have $\psi\in\Psi[T_1](\RR,\RR)$.
By part (iv) of Theorem \ref{Thm_M_est_uniq}, for each $x,y\in \RR$ with $\vartheta_1(x)<\vartheta_1(y)$,
 the function \eqref{function_newhanyados0} is increasing.
Since $\psi\in\Psi[T_1](\RR,\RR)$, by part (a) of the present proposition, we have $f(z)>0$ for each $z>0$ and $\vartheta_1(x)=x$.
Consequently, for each $x<y$, $x,y\in\RR$, the function (given by \eqref{function_newhanyados0})
 \begin{align}\label{function_newhanyados00}
   (x,y)\ni t \mapsto  -\frac{\psi(x,t)}{\psi(y,t)}
 \end{align}
 is increasing.
Hence the statement of part (b) follows by the following observation (that we check below):
 provided that $f(z)>0$ for each $z>0$,
 the function \eqref{function_newhanyados00} is (strictly) increasing
 for each $x<y$, $x,y\in\RR$ if and only if $f$ is (strictly) increasing.
Since the function $\psi$ given in \eqref{help_psi_mathieu}
 depends only on $x-t$, it is enough to check that
 \begin{align}\label{help_psi_mathieu_2}
   \text{the function} \quad(0,z)\ni t\mapsto \frac{\psi(0,t)}{\psi(z,t)} = -\frac{f(t)}{f(z-t)}
    \quad \text{is (strictly) decreasing for each $z>0$}
 \end{align}
 holds if and only if $f$ is (strictly) increasing.
Indeed, for each $x<y$, $x,y\in\RR$, and $t\in(x,y)$, we have
 \[
   \frac{\psi(x,t)}{\psi(y,t)}
     = -\frac{f(t-x)}{f(y-t)}
     = -\frac{f(t-x)}{f(y-x-(t-x))}
     = \frac{\psi(0,t-x)}{\psi(y-x,t-x)},
     \qquad t-x\in(0,y-x).
 \]
Thus, the property \eqref{help_psi_mathieu_2} holds if and only if
  \[
   \frac{f(s)}{f(z-s)} \, (<)\leq \, \frac{f(t)}{f(z-t)}\qquad  \text{for each $s,t,z\in\RR$ with $0<s<t<z$,}
  \]
  which is equivalent to
  \begin{align}\label{help18}
   f(s)f(z-t) \, (<)\leq \,  f(t)f(z-s)  \qquad \text{for each $s,t,z\in\RR$ with $0<s<t<z$.}
  \end{align}
Using the nonnegativity of $f$, it yields that \eqref{help_psi_mathieu_2} holds if and only if
  \begin{align}\label{help19}
   f(s) \, (<)\leq \, f(t) \qquad \text{for each $s,t\in\RR$ with $0<s<t$,}
  \end{align}
  i.e., $f$ is (strictly) increasing on $\RR_+$.
Indeed, if \eqref{help_psi_mathieu_2} holds, then \eqref{help18} holds as well, and, by choosing $z=s+t$, we get that $f(s)^2(<)\leq f(t)^2$,
 which implies \eqref{help19}, since $f$ is nonnegative.
If \eqref{help19} holds, then for each $s,t,z\in\RR$ with $0<s<t<z$, we have $0\leq f(s) \, (<)\leq \, f(t)$
 and $0\leq f(z-t) \, (<)\leq \, f(z-s)$, implying \eqref{help18}, and hence \eqref{help_psi_mathieu_2} as well.

Parts (c) and (d):
Let us suppose that $\psi\in\Psi[T_2^{\blambda}](\RR,\RR)$ for all $\blambda\in\Lambda_2$.
In particular, $\psi\in\Psi[T_1](\RR,\RR)$, and hence, by part (a) of the present proposition,
 we have $f(z)>0$ for each $z>0$ and $\vartheta_1(x)=x$, $x\in\RR$.
Consequently, by part (v) of Theorem \ref{Thm_M_est_uniq}, for each $x<y$, $x,y\in\RR$, the function
 \eqref{function_newhanyados00} is strictly increasing.
By the proof of part (b) of the present proposition, it implies that $f$ is strictly increasing, as desired.
We give an alternative proof as well. Let $s,t\in\RR_+$ with $s<t$, $x_1\in\RR$, $x_2:=x_1+s+t$, and $r:=x_1+s$.
Then $\frac{1}{2}(x_1+x_2)=x_1+\frac{s+t}{2}>x_1+s=r$, yielding that
 \[
   \psi\Big(x_1,\frac{x_1+x_2}{2}\Big) + \psi\Big(x_2,\frac{x_1+x_2}{2}\Big)
      = - f\Big(\frac{s+t}{2}\Big) + f\Big(\frac{s+t}{2}\Big) =0.
 \]
Hence $\vartheta_{2,\psi}(x_1,x_2)=\frac{1}{2}(x_1+x_2)$, and, since $r<\frac{1}{2}(x_1+x_2)$, we have that
 $\psi(x_1,r)+\psi(x_2,r)>0$.
Since $\psi(x_1,r) + \psi(x_2,r) = -f(s)+f(t)$, we get that $f(s)<f(t)$, as desired.
Conversely, let us suppose that $f$ is strictly increasing.
Then for each $x\in\RR$, the function \ $\RR\ni t\mapsto \psi(x,t)$ \ is strictly decreasing.
Hence Proposition \ref{Pro_M_est_uniq} implies that $\psi\in\Psi[T_n^\blambda](\RR,\RR)$
 for each $n\in\NN$ and $\blambda\in\Lambda_n$ (in particular, $\psi\in\Psi[T_2^\blambda](\RR,\RR)$-function
 for each $\blambda\in\Lambda_n$), as desired.
\proofend

\section*{Acknowledgements}

We acknowledge the valuable suggestions from the two referees.

\small
\bibliographystyle{plain}
\bibliography{psi_estimator_bib}

\begin{thebibliography}{10}

\bibitem{AliTib}
A.~Ali and R.~J. Tibshirani.
\newblock The generalized lasso problem and uniqueness.
\newblock {\em Electron. J. Stat.}, 13(2):2307--2347, 2019.

\bibitem{BelKlaMul}
F.~Bellini, B.~Klar, and A.~M\"{u}ller.
\newblock Expectiles, {O}mega ratios and stochastic ordering.
\newblock {\em Methodol. Comput. Appl. Probab.}, 20(3):855--873, 2018.

\bibitem{Cat12}
O.~Catoni.
\newblock Challenging the empirical mean and empirical variance: a deviation
  study.
\newblock {\em Ann. Inst. Henri Poincar\'{e} Probab. Stat.}, 48(4):1148--1185,
  2012.

\bibitem{CheJinLiXu}
P.~Chen, X.~Jin, X.~Li, and L.~Xu.
\newblock A generalized {C}atoni's {${\rm M}$}-estimator under finite
  {$\alpha$}-th moment assumption with {$\alpha\in(1,2)$}.
\newblock {\em Electron. J. Stat.}, 15(2):5523--5544, 2021.

\bibitem{Cla}
B.~R. Clarke.
\newblock Uniqueness and {F}r\'{e}chet differentiability of functional
  solutions to maximum likelihood type equations.
\newblock {\em Ann. Statist.}, 11(4):1196--1205, 1983.

\bibitem{Cla2}
B.~R. Clarke.
\newblock The selection functional.
\newblock {\em Probab. Math. Statist.}, 11(2):149--156 (1991), 1990.

\bibitem{DimFisZie}
T.~Dimitriadis, T.~Fissler, and J.~Ziegel.
\newblock Characterizing {$M$}-estimators.
\newblock {\em Biometrika}, 111(1):339--346, 2024.

\bibitem{GasPap}
L.~Gasi\'nski and N.~S. Papageorgiou.
\newblock {\em Exercises in {A}nalysis. {P}art 2. {N}onlinear {A}nalysis}.
\newblock Problem Books in Mathematics. Springer, Cham, 2016.

\bibitem{GruPal}
R.~Gr\"{u}nwald and {\relax Zs}.~P\'{a}les.
\newblock On the equality problem of generalized {B}ajraktarevi\'{c} means.
\newblock {\em Aequationes Math.}, 94(4):651--677, 2020.

\bibitem{HamHenRon}
F.~Hampel, C.~Hennig, and E.~Ronchetti.
\newblock A smoothing principle for the {H}uber and other location
  {$M$}-estimators.
\newblock {\em Comput. Statist. Data Anal.}, 55(1):324--337, 2011.

\bibitem{Hub64}
P.~J. Huber.
\newblock Robust estimation of a location parameter.
\newblock {\em Ann. Math. Statist.}, 35:73--101, 1964.

\bibitem{Hub67}
P.~J. Huber.
\newblock The behavior of maximum likelihood estimates under nonstandard
  conditions.
\newblock In {\em Proc. {F}ifth {B}erkeley {S}ympos. {M}ath. {S}tatist. and
  {P}robability ({B}erkeley, {C}alif., 1965/66), {V}ol. {I}: {S}tatistics},
  pages 221--233. Univ. California Press, Berkeley, Calif., 1967.

\bibitem{JanMer}
S.~Jankovi\'{c} and M.~Merkle.
\newblock A mean value theorem for systems of integrals.
\newblock {\em J. Math. Anal. Appl.}, 342(1):334--339, 2008.

\bibitem{KoeBas}
R.~Koenker and Jr.~G. Bassett.
\newblock Regression quantiles.
\newblock {\em Econometrica}, 46(1):33--50, 1978.

\bibitem{Kos}
M.~R. Kosorok.
\newblock {\em Introduction to Empirical Processes and Semiparametric
  Inference}.
\newblock Springer Series in Statistics. Springer, New York, 2008.

\bibitem{KraZah}
V.~Kr\"atschmer and H.~Z\"ahle.
\newblock Statistical inference for expectile-based risk measures.
\newblock {\em Scand. J. Stat.}, 44(2):425--454, 2017.

\bibitem{Lan}
K.~Lange.
\newblock {\em Numerical Analysis for Statisticians}.
\newblock Statistics and Computing. Springer, New York, second edition, 2010.

\bibitem{Mat}
T.~Mathieu.
\newblock Concentration study of {M}-estimators using the influence function.
\newblock {\em Electron. J. Stat.}, 16(1):3695--3750, 2022.

\bibitem{NewPow}
W.~K. Newey and J.~L. Powell.
\newblock Asymmetric least squares estimation and testing.
\newblock {\em Econometrica}, 55(4):819--847, 1987.

\bibitem{Pal2003}
{\relax Zs}.~P\'{a}les.
\newblock On approximately convex functions.
\newblock {\em Proc. Amer. Math. Soc.}, 131(1):243--252, 2003.

\bibitem{PasRei}
R.~Passeggeri and N.~Reid.
\newblock A universal robustification procedure.
\newblock Available on {\it arXiv}: 2206.06998, 2022.

\bibitem{Rey}
W.~J.~J. Rey.
\newblock {\em Introduction to Robust and Quasirobust Statistical Methods}.
\newblock Universitext. Springer-Verlag, Berlin, 1983.

\bibitem{ShaDenRus}
A.~Shapiro, D.~Dentcheva, and A.~Ruszczy\'nski.
\newblock {\em Lectures on Stochastic Programming}, volume~9 of {\em MPS/SIAM
  Series on Optimization}.
\newblock Society for Industrial and Applied Mathematics (SIAM), Philadelphia,
  PA; Mathematical Programming Society (MPS), Philadelphia, PA, 2009.

\bibitem{Tib}
R.~J. Tibshirani.
\newblock The lasso problem and uniqueness.
\newblock {\em Electron. J. Stat.}, 7:1456--1490, 2013.

\bibitem{Vaa}
A.~W. van~der Vaart.
\newblock {\em Asymptotic {S}tatistics}, volume~3 of {\em Cambridge Series in
  Statistical and Probabilistic Mathematics}.
\newblock Cambridge University Press, Cambridge, 1998.

\end{thebibliography}

\end{document}